\documentclass[preprint]{elsarticle}

\usepackage{hyperref}
\usepackage[centertags]{amsmath}
\usepackage{amsfonts}
\usepackage{amssymb}
\usepackage{amsmath, bm}
\usepackage[ruled]{algorithm2e}
\usepackage{enumitem}
\usepackage{multicol}
\usepackage{multirow}
\usepackage{makecell}
\usepackage{ragged2e}
\usepackage{color}
\usepackage{subfigure}

\usepackage{amsthm}

\newtheorem{remark}{Remark}[section]

\journal{Journal}
\SetKwBlock{pre}{Pre-training step:}{end}
\SetKwBlock{formal}{Formal training step:}{end}
\numberwithin{table}{section}

\numberwithin{equation}{section}
\numberwithin{figure}{section}
\begin{document}
\begin{frontmatter}
\title{
A deep learning method for multi-material  diffusion problems based on physics-informed neural networks
}
\author[1]{Yanzhong Yao}
\author[1,2]{Jiawei Guo\corref{cor1}}
\ead{guojiawei20@gscaep.ac.cn}
\author[1]{Tongxiang Gu}

\cortext[cor1]{Corresponding author}
\address[1]{Laboratory of Computational Physics, Institute of Applied Physics and
Computational Mathematics, Beijing 100088, China}
\address[2]{Graduate School of China Academy of Engineering Physics, Beijing 100088, China}

\begin{abstract}
Given the facts of the extensiveness of multi-material diffusion problems and the inability of the standard PINN(Physics-Informed Neural Networks) method for such problems, in this paper we present a novel PINN method that can accurately solve the multi-material  diffusion equation. The new method applies continuity conditions at the material interface derived from the property of the diffusion equation, and combines the distinctive spatial separation strategy and the loss term normalization strategy to solve the problem that the residual points cannot be arranged at the material interface,
the problem that it is difficult to express non-smooth functions with a single neural network,
and the problem that the neural network is difficult to optimize the loss function with different magnitudes of loss terms,
which finally provides the available prediction function for a class of multi-material diffusion problems. 
Numerical experiments verify the robustness and effectiveness of the new method.

\end{abstract}

\begin{keyword}
multi-material diffusion equation
\sep deep learning method
\sep physics-informed neural networks
\sep flux continuity condition
\sep domain separation strategy.
\end{keyword}

\end{frontmatter}


\section{Introduction}
Diffusion equations are an important class of partial differential equations that need to be studied in many applications such as groundwater seepage \cite{Guo}, oil reservoir simulation \cite{ILLARIONOV}, nuclear reactions \cite{Lindl}, etc. They can formulate both the diffusion process of material concentration in space, and 
the energy transfer process in radiative heat conduction problems \cite{Yao2,Mishra}. Their numerical solution methods have been a hot research topic in the field of scientific and engineering computation.

In recent years, the technology of deep learning for solving partial differential equations has developed rapidly. Researchers have developed a new method for solving partial differential equations by introducing Physical Information into Neural Networks, which is known as the PINN method \cite{Karniadakis,RAISSI1,Raissi4,Raissi_sci}. The (continuous) PINN method described in Ref.~\cite{RAISSI1} will be referred to as {\itshape{the standard PINN method}} in the remainder of this paper. The PINN method takes the definite solution condition as the supervised learning part, and combines the automatic differential technology to take the approximation degree of the governing equation as the residual part. Two parts together form the loss function as the training objective of the network prediction. This technique makes the network output obtained by the optimization algorithm not only satisfy the definite solution condition, but also satisfy the governing equation.

Compared with the traditional numerical methods, including the finite element method (FEM) and the finite volume method, the PINN method (FVM) has advantages in some aspects, such as unnecessary grid generation, adaptation to high-dimensional problems, and so on.
However, the PINN method still faces many challenges in practical applications, one of which is how to efficiently solve heterogeneous diffusion equations on the computational domain involving multiple materials.

For the multi-material diffusion problem, there will be material interfaces, and because the materials on both sides of the interface have different physical properties, their diffusion coefficients or heat conduction coefficients will have large differences, and they will have jumps at the interface, and such jumps will cause the derivatives of the solution function to necessarily have jumps at the interface, i.e., the solution function is not continuously differentiable at the interface, and the second-order partial derivatives of the solution function will not exist at all. This problem poses two major difficulties for the standard PINN method: \emph{(1) the PINN method generally produces a smooth prediction function, so it is difficult for the standard PINN method to obtain a prediction function that is not continuously differentiable \cite{Lu};  (2) since the solution function does not have second-order partial derivatives at the interface, the standard PINN method based on the automatic differentiation technique cannot incorporate the sampled points at the interface as residual points in the training of a single neural network \cite{HE1}. This inevitably leads to an inaccurate solution near the interface.}
However, it is well known that the structure of the solution near the interface is generally very complex, and the computational accuracy for this position can seriously affect the overall computational accuracy of the computational model. Obtaining high accuracy numerical solutions near the interface is a very challenging but necessary task for any numerical method.

To solve the equation with a non-smooth solution at interfaces using the PINN method, the following two strategies are most commonly used.
One is to accept the fact that the neural network will make an incorrect prediction near the interface. To get useful predictions in the part away from the interface, the points at the interface are not sampled, or their residuals are weighted so that the loss near the interface tends to 0. 
In Ref.~\cite{Xie}, Xie et al. designed a weighted function to handle the possible jump of the diffusion coefficient across the material interface. 
The other idea is to use different neural networks for each subdomain, so that the outputs of the multiple networks can express the functions with non-smoothness at the interface.
In Ref.~\cite{HE1}, He et al. pointed out that it is inefficient to use only one neural network structure to capture the interface, and they proposed a piece-wise neural network structure and an adaptive resampling strategy to solve the elliptic interface problems. As a result, the PINN method in combination with domain decomposition techniques has received increasing attention.
For solving nonlinear PDEs in complex geometry domains, extended PINNs (XPINNs) based on domain decomposition techniques have been proposed in Ref.~\cite{Jagtap2}.
In Ref.~\cite{Dwivedi}, a Distributor-PINN method (DPINN) was proposed, which decomposes the computational domain into several disjoint subdomains and trains several sub-networks simultaneously by minimizing the loss function, and each sub-network is used to approximate the solution on a subdomain.
Deep DDM in \cite{LiWy} is another PINN method based on domain decomposition techniques to solve two-dimensional elliptic interface problems.
In Ref.\cite{Wu}, the author showed that the treatments in the above domain decomposition-based PINN methods could also suffer from convergence issues or have the drawback of low accuracy by the experiments, and they developed a dynamic weighting strategy that adaptively assigns appropriate weights to different terms in the loss function based on the multi-gradient descent algorithm \cite{Dee}. 

In view of the facts of the extensiveness of multi-material diffusion problems and the inability of the standard PINN for such problems, this paper attempts to propose effective strategies to overcome the above two difficulties, including that a single neural network is difficult to express the function with different derivatives on two sides of the interface, and the standard PINN method cannot arrange effective sampling points on the interface, so as to construct a novel PINN method named as \emph{DS-PINN}, which can accurately solve the multi-material diffusion equation using a single neural network. In addition, we develop a normalization strategy for the loss terms and present the \emph{nDS-PINN} method, which further improves the prediction accuracy.

The rest of this paper is organized as follows: in section 2, we do some preliminary work by giving the governing equation of the multi-material  diffusion problem and its standard PINN form; in section 3, we discuss in detail how to improve the standard PINN to obtain an available prediction function for the multi-material diffusion problem; in section 4, we give several numerical examples to verify the effectiveness of our method; the conclusion about the new PINN method is drawn in the last section.

\section{Preliminaries} 

In this section, we first describe a class of multi-material diffusion problems and then give the standard PINN method for solving them.

\subsection{Physical model and governing equation}

A class of linear diffusion problems on the domain containing different materials can be formulated as follows:


\begin{equation}\label{2.1}
-\nabla\cdot \kappa(X)\nabla u = Q\left(X\right), X\in\Omega,
\end{equation}
with the Dirichlet boundary condition
\begin{equation}\label{bdc}
u(X)=g(X),X\in\partial\Omega,
\end{equation}
where $u=u(X)$ is the function to be solved, 
$\Omega$ is an open domain in $\mathbb{R}^d$ with the boundary $\partial \Omega$.
For multi-material problems, $\Omega$ is composed of several subdomains containing different materials.
The source term $Q(X)$ and the boundary condition $g(X)$ are bounded in their domains of definition. 
The boundary condition can also be of Robin or mixed type. 
Note that, the diffusion coefficient
\begin{equation}\label{mmk}
\kappa(X)=\kappa_i(X),~\text{for}~X\in\Omega_i,~i=1,2,\cdots,N,
\end{equation}
where $\Omega_i$ denotes a subdomain containing a certain material, 
$\kappa_i(X)$ is the diffusion coefficient on that subdomain,
and $N$ denotes the number of material type. 
$\kappa_i(X),i=1,\dots,N$ are smooth functions on their respective subdomains, but they may not be equal at the material interface.

To simplify the description, this paper mainly discusses two-dimensional problems, and the same idea can also be applied to three-dimensional problems.

\begin{figure}[h]
\centering
\includegraphics[width=0.35\textwidth]{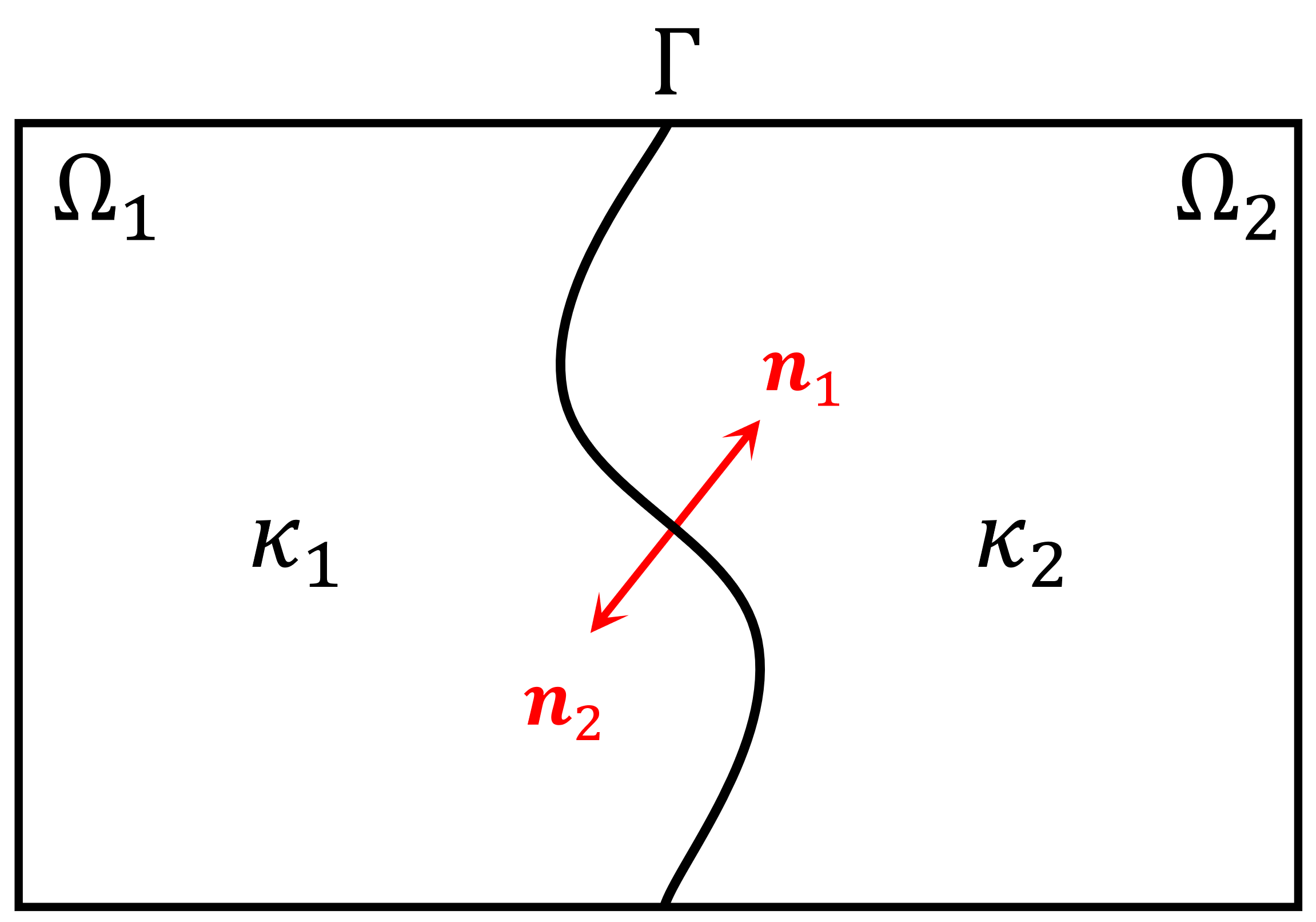}
\caption{A 2D domain with two materials.}
\label{Fig1.1}
\end{figure}

Figure \ref{Fig1.1} shows a 2D domain containing two materials, and $\Gamma$ is the material interface. If $\kappa_1(X)\neq\kappa_2(X), X\in \Gamma $, then according to Eq.~\eqref{2.1}, one can get 
\begin{equation}
     \nabla u(X) | _{\Gamma^-}
    \neq\nabla u(X) | _{\Gamma^+}.
\end{equation}
Thus, the solution $u$ is not continuously differentiable at the material interface $\Gamma$,
and its second-order partial derivatives
are absent.

\subsection{Standard PINN method for diffusion equations}

The PINN method produces the prediction function $u_\theta$ as an approximation of the solution of Eqs.~\eqref{2.1}-\eqref{mmk}, where $\theta$ denotes the neural network parameters, including the weights and biases of the neural networks.
Later in this paper, we refer to $u_\theta$ as \emph{the prediction}.

For the standard PINN method, $\theta$ are obtained by optimising the loss function, and the loss function consists of two parts as follows:
\begin{equation}\label{loss_std}
\mathcal{L}(\theta;\Sigma)=w_b\mathcal{L}_b(\theta;\tau_b)+w_r\mathcal{L}_r(\theta;\tau_r),
\end{equation}
where $w_b$ and $w_r$ are the weights for the two parts of the loss function, respectively,
and the \emph{supervised loss term} and the \emph{residual loss term} are
\begin{align}
&\mathcal{L}_{b}(\theta;\tau_{b}) = \frac{1}{N_b}\sum_{i=1}^{N_b}\left|u_\theta(X_i)-g(X_i)\right|^2,\\
&\mathcal{L}_{r}(\theta;\tau_{r})=\frac{1}{N_r}\sum_{i=1}^{N_{r}}\left|-\nabla\cdot \kappa(X_i)\nabla u(X_i) - Q\left(X_i\right)
\right|^2.  \label{Lr_loss}
\end{align}
$\Sigma = \{\tau_b,\tau_r\}$ denotes the training data set, where
$\tau_b=\{\left(X_i,g(X_i)\right)|X_i\in\partial\Omega\}_{i=1}^{N_b}$ is the labeled data set,  
and $\tau_{r}=\{X_i \in\Omega\}_{i=1}^{N_{r}}$ is the residual data set.
$N_b$ and $N_{r}$ denote the number of boundary sampling points on $\partial\Omega$ and the number of inner sampling points in $\Omega$, respectively.

Suppose 
\begin{equation}\label{MIN}
\bar{\theta} = \arg\min_\theta\mathcal{L}(\theta;\Sigma),
\end{equation}
which can be obtained by some optimization methods, and then $u_{\bar{\theta}}$ is the approximation of the unknown function $u$.

The partial differential operators, such as $u_x$ and $u_{xx}$, can be implemented using automatic differentiation (AD). This can be easily realized in the deep learning framework like the PyTorch \cite{PyTorch} or Tensorflow \cite{TF}.
\begin{remark}
For unsteady diffusion problems 
 \begin{equation}\label{unsteadyE}
\begin{cases}
    u_t-\nabla\cdot \kappa(t,X)\nabla u = Q\left(t,X\right), t\in(0,T],&~X\in\Omega,\\
u(t,X)=g(t,X),t\in(0,T],&~X\in\partial\Omega,\\
u(0,X)=\phi(X),&~X\in\Omega,
\end{cases}
\end{equation}
the method to be discussed in this paper is also applicable. At this case, the loss function requires adding a loss term to reflect the degree of approximation of the initial condition.
\end{remark}

\section{An improved PINN method for solving multi-material diffusion problems}

In this section, we investigate a deep learning method for solving the multi-material diffusion equations \eqref{2.1}-\eqref{mmk},  and present an improved PINN  using interface connection conditions and the domain separation strategy, which is called DS-PINN. At the end of this section, we further improve the performance of the DS-PINN method by introducing the normalization strategy, which is denoted as nD-PINN.

It is well known that the standard PINN method, under the assumption that the solution of the equation is sufficiently smooth, uses automatic differentiation techniques to compute the residual loss term. Since the diffusion coefficients of Eq.~\eqref{2.1} are discontinuous at the material interface, $u$ is not continuously differentiable  and its second-order derivatives do not exist at the interface, which means that one cannot sample residual points at the interface, and thus the equation information at the interface is lost.
If this issue is not properly addressed, the prediction obtained from the standard PINN will have a large uncertainty at the material interface, resulting in an unreliable result.
The following subsections provide several strategies for dealing with this problem.


\subsection{Introducing material interface continuity conditions into the standard PINN}\label{linkcond}

Since the second derivative of the function $u$ does not exist at the material interface, the residual error of the sampling point at the interface cannot be calculated according to Eq.~\eqref{Lr_loss}. To compensate for this deficiency, we can add new loss terms to the loss function according to the properties that Eq.~\eqref{2.1} satisfies at the interface, so that the prediction function of the neural network can reasonably reflect the behavior of the solution at the interface.


According to the property of the diffusion equation, the following two conditions should be satisfied at the interface $\Gamma$:
\begin{align} 
    \left.u(X)\right|_{\Gamma^-}&= \left.u(X)\right|_{\Gamma^+},\label{cu}\\
    \left.-\kappa_1(X)\nabla u(X)\right|_{\Gamma^-}\cdot\boldsymbol{n}_1&=-\left(\left.-\kappa_2(X)\nabla u(X)\right|_{\Gamma^+}\cdot\boldsymbol{n}_2\right), \label{cFlux}
\end{align}
where $\left.\cdot\right|_{\Gamma^-}$ and $\left.\cdot\right|_{\Gamma^+}$  represent the corresponding function values of approaching any point $X$ on the interface $\Gamma$ from $\Omega_1$ and $\Omega_2$, respectively.
$\boldsymbol{n}_1$ and $\boldsymbol{n}_2$ are the outer normal directions of the corresponding subdomain, as shown in the Figure \ref{Fig1.1}. 


Eqs.~\eqref{cu}-\eqref{cFlux} are  called the continuity conditions, which include the solution continuity and  the flux continuity. 


Define $[\![\mathcal F(X)]\!]_{\Gamma}:=\mathcal F(X)|_{\Gamma^+} - \mathcal F(X)|_{\Gamma^-}$ to denote the jump of $\mathcal F(X)$ across the material interface. Then the continuity conditions ~\eqref{cu}-\eqref{cFlux} can be rewritten as follows:
\begin{align} 
    [\![u(X)]\!]_{\Gamma}&=0,\label{cu_new}\\
    [\![-\kappa(X)\nabla u(X)\cdot\boldsymbol{n} ]\!]_{\Gamma}&=0.\label{cFlux_new}
\end{align}

In fact, the continuity conditions of the material interface are the hypothetical conditions for the derivation of Eq.~\eqref{2.1}. 
However, these two conditions can also be obtained from Eq.~\eqref{2.1}. 
The solution continuity condition \eqref{cu} is the assumed condition, which is a necessary condition for the governing equation \eqref{2.1} to hold. 
Next, we give the derivation of the flux continuity condition \eqref{cFlux}.

\begin{figure}
\centering
\includegraphics[width=0.4\textwidth]{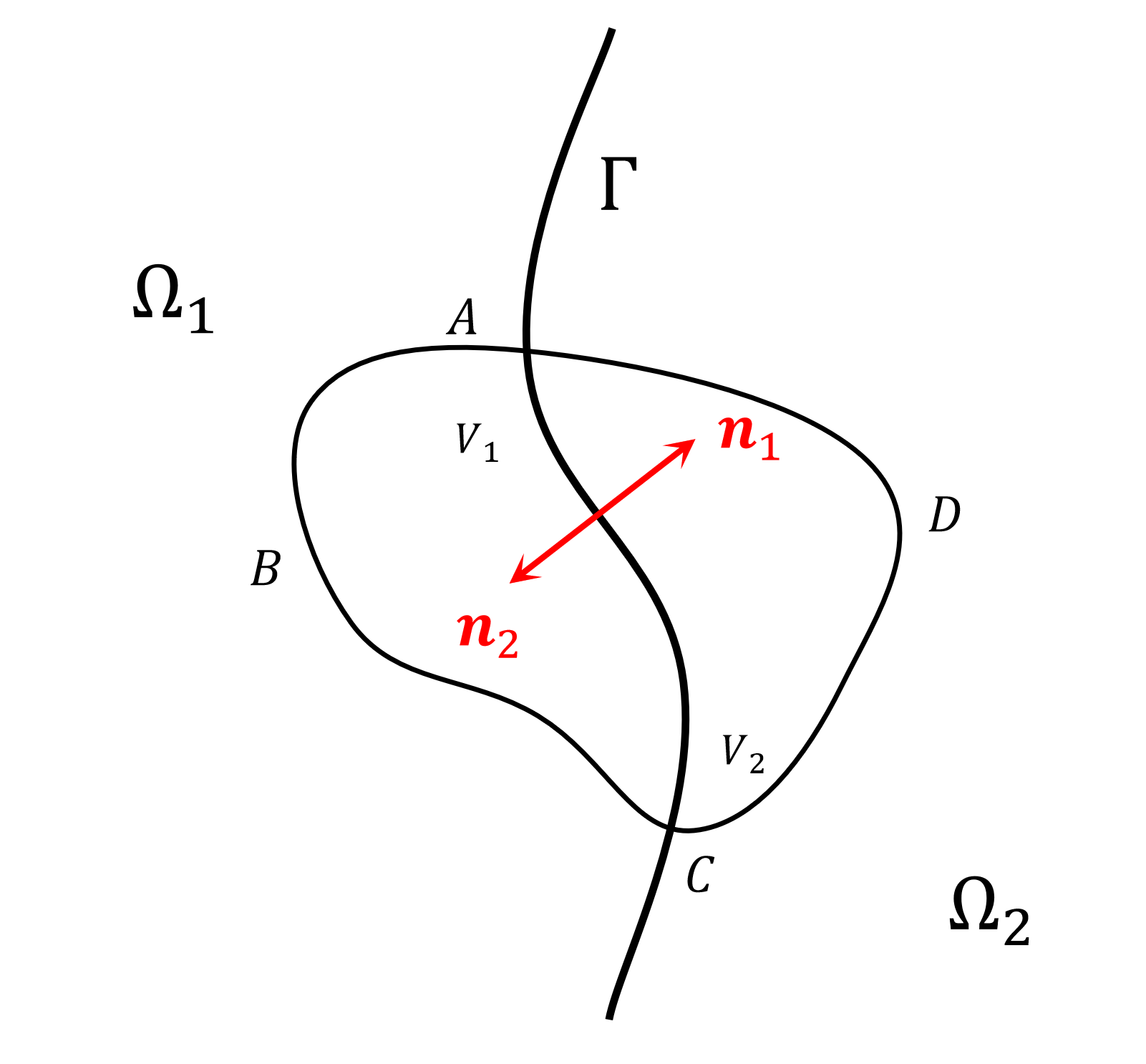}
\caption{A control volume $V(V_1\cup V_2)$ containing part of the interface $\Gamma$.}
\label{Fig_flux}
\end{figure}

Suppose that $V$ is an arbitrary control volume containing part of the interface $\Gamma$ which divides it into two parts $V_1$ and $V_2$, as shown in Figure \ref{Fig_flux}.

Integrating Eq. \eqref{2.1}  over the volume $V_1$, we obtain
\begin{equation}\label{int1}
\int_{V_1}\left(-\nabla\cdot\kappa_1(X) \nabla u\right)\mathrm{d}V = \int_{V_1} Q(X) \mathrm{d}V.
\end{equation}
According to the divergence theorem, Eq.~\eqref{int1} can be rewritten as follows:
\begin{align}
\int_{V_1}\left(-\nabla\cdot\kappa_1(X) \nabla u\right)\mathrm{d}V = \oint_{\partial V_1}&\left(-\kappa_1(X) \nabla u\right)\cdot\boldsymbol{n}_1\mathrm{d}S\notag\\
=\int_{\overrightarrow{ABC} + \overrightarrow{CA}}&\left(-\kappa_1(X) \nabla u\right)\cdot\boldsymbol{n}_1\mathrm{d}S=\int_{V_1} Q(X) \mathrm{d}V.\label{flux1}
\end{align}
Analogously, we obtain the similar result for $V_2$
\begin{align}
\int_{\overrightarrow{CDA} + \overrightarrow{AC}}\left(-\kappa_2(X)\nabla u\right)\cdot\boldsymbol{n}_2\mathrm{d}S = \int_{V_2} Q(X) \mathrm{d}V .\label{flux2}
\end{align}
Integrating Eq.~\eqref{2.1} over the entire control volume $V$ and using the divergence theorem, we have
\begin{align}
\int_{\overrightarrow{ABC} + \overrightarrow{CDA}}\left(-\kappa(X)\nabla u\right)\cdot\boldsymbol{n}\mathrm{d}S = \int_{V} Q(X) \mathrm{d}V.\label{fluxall}
\end{align}
By combining Eqs.~\eqref{flux1},~\eqref{flux2} and \eqref{fluxall}, we get the flux continuity formula 
\begin{equation}\label{flux_lcont}
\int_{\overrightarrow{CA}}\left(-\kappa_1(X) \nabla u\right)\cdot\boldsymbol{n}_1\mathrm{d}S + \int_{\overrightarrow{AC}}\left(-\kappa_2(X) \nabla u\right)\cdot\boldsymbol{n}_2\mathrm{d}S = 0.
\end{equation}
Given the arbitrariness of $V$, we can obtain the flux continuity condition \eqref{cFlux}.

Eq.\ref{flux_lcont} is used in many papers constructing discrete schemes for heterogeneous diffusion equations, such as \cite{Aavatsmark, Miao2, Brezzi}. 

\begin{remark}
    The equations should satisfy the continuity condition at the interface, so we want to reflect this property of the predicted solution near the interface by adding new loss terms to the loss function. However, these two conditions cannot be directly applied to the training of the PINN. The reason is that the prediction function obtained from a single PINN training is continuously differentiable, i.e., there is only one unique derivative at each position, so the solution continuity condition is naturally satisfied, while the flow continuity condition cannot be achieved.
\end{remark}



\subsection{Applying domain separation strategy to compute derivatives on both sides of the material interface}\label{domsep}



To characterize the derivative discontinuity property at the interface, a natural idea is to train this model using two sets of neural networks linked by the interface connection conditions. However, this strategy faces some difficulties: \emph{(1) it requires multiple sets of networks for the multi-material model; (2) it is more difficult to design and implement optimization algorithms; and (3) it generally requires iteration between neural networks and the convergence is difficult to guarantee.}


Under the premise of using only one neural network to obtain different derivative values on both sides of the interface, an intuitive idea is to separate the two domains divided by the interface by a certain distance, and the material interface becomes the boundaries of two sub-domains with a certain interval, so that each point on the interface becomes two points belonging to different locations in space, so it is logical that they can have different derivative values.

\begin{figure}[h]
\centering
\subfigure[]{\includegraphics[width=0.4\textwidth]{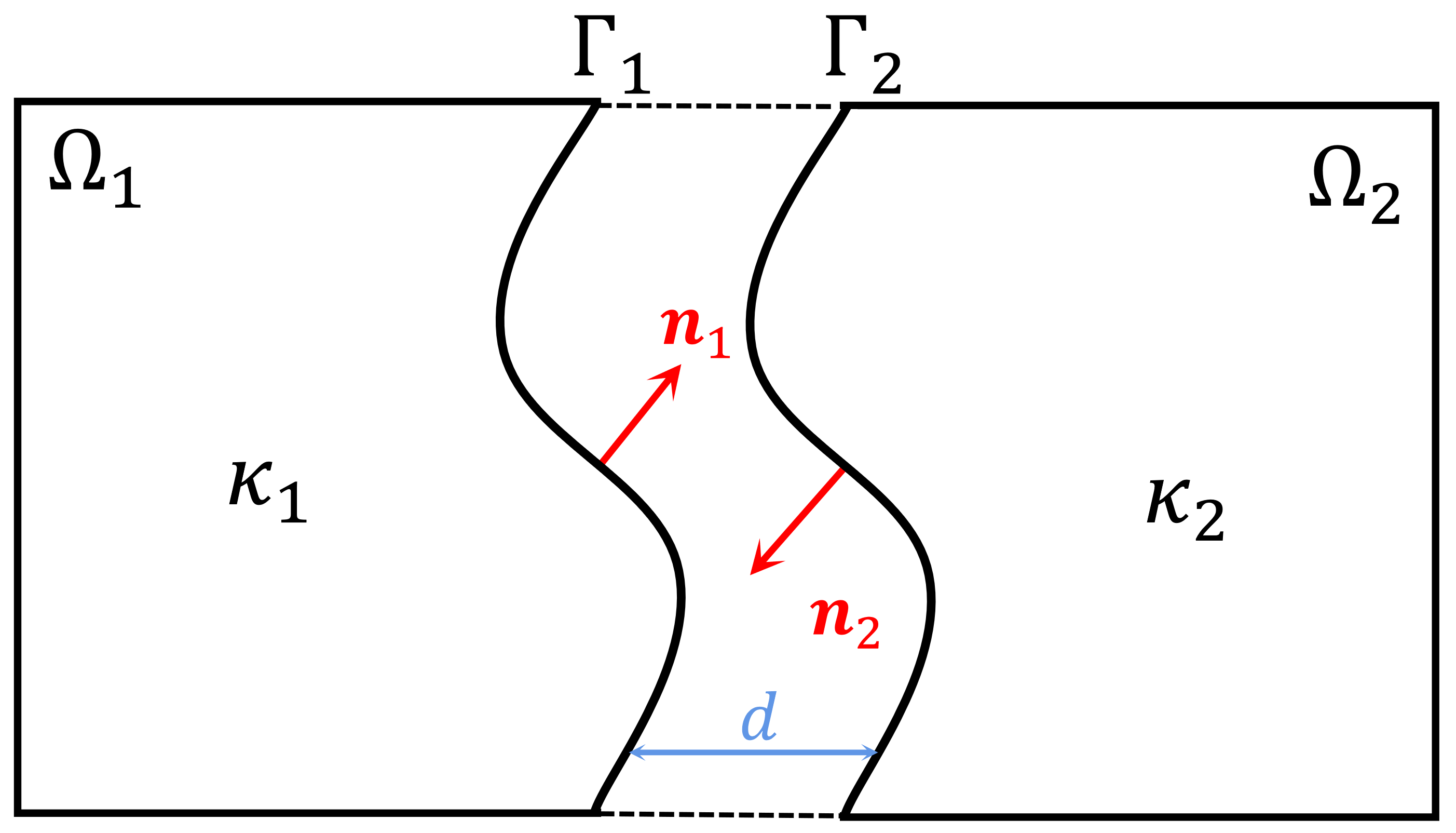}}
\subfigure[]{\includegraphics[width=0.415\textwidth]{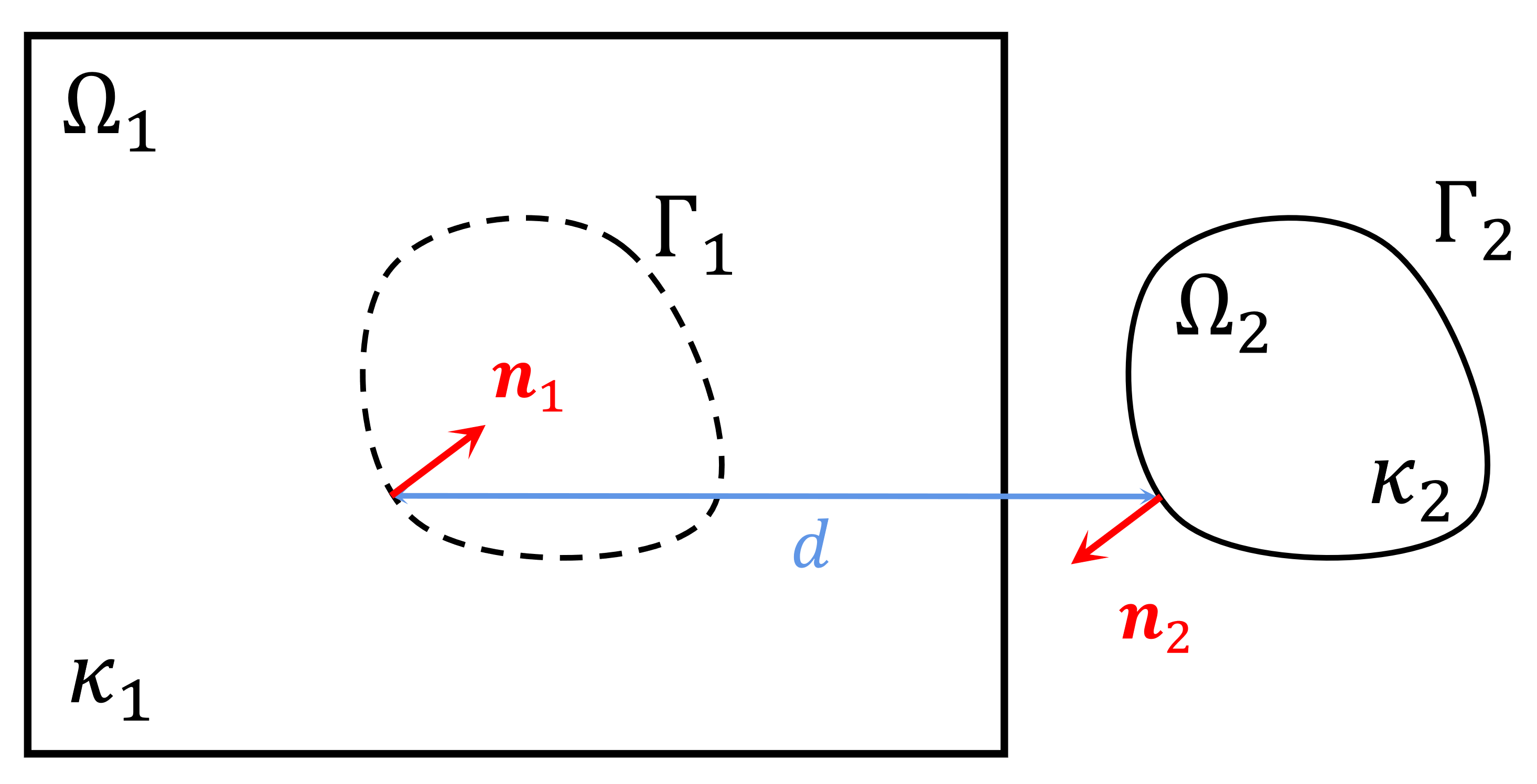}}\\
\caption{Two ways of separating domains. }
\label{Fig2.1}
\end{figure}

Figure \ref{Fig2.1} shows the respective domain separation strategies for the two types of material interfaces.
It is easy to generalize this strategy to the case of multiple materials.

With regard to the domain separation strategy, the following points should be noted:

\begin{itemize}
    \item The subdomains cannot overlap because only one neural network is used.
    
 \item There are no strict limitations on the distance and direction of the separation. 
 If the difference of diffusion coefficients between the two sides of the interface is large, then we should choose a larger distance. 
 In Sect.~\ref{NEX}, the performance is tested with different separation distances ${d}$.

 \item After implementing the separation strategy, the material interface changes from $\Gamma$ to two spatially separated boundaries $\Gamma_1$ and $\Gamma_2$, on which the interface continuity conditions \eqref{cu},\eqref{cFlux} are imposed.

 \item The training points are sampled based on all subdomains after separation, and the sampling points on the material interface should be consistent between the matching subdomains to impose the interface continuity condition.
\end{itemize}

\begin{figure}[h]
\centering
\includegraphics[width=0.5\textwidth]{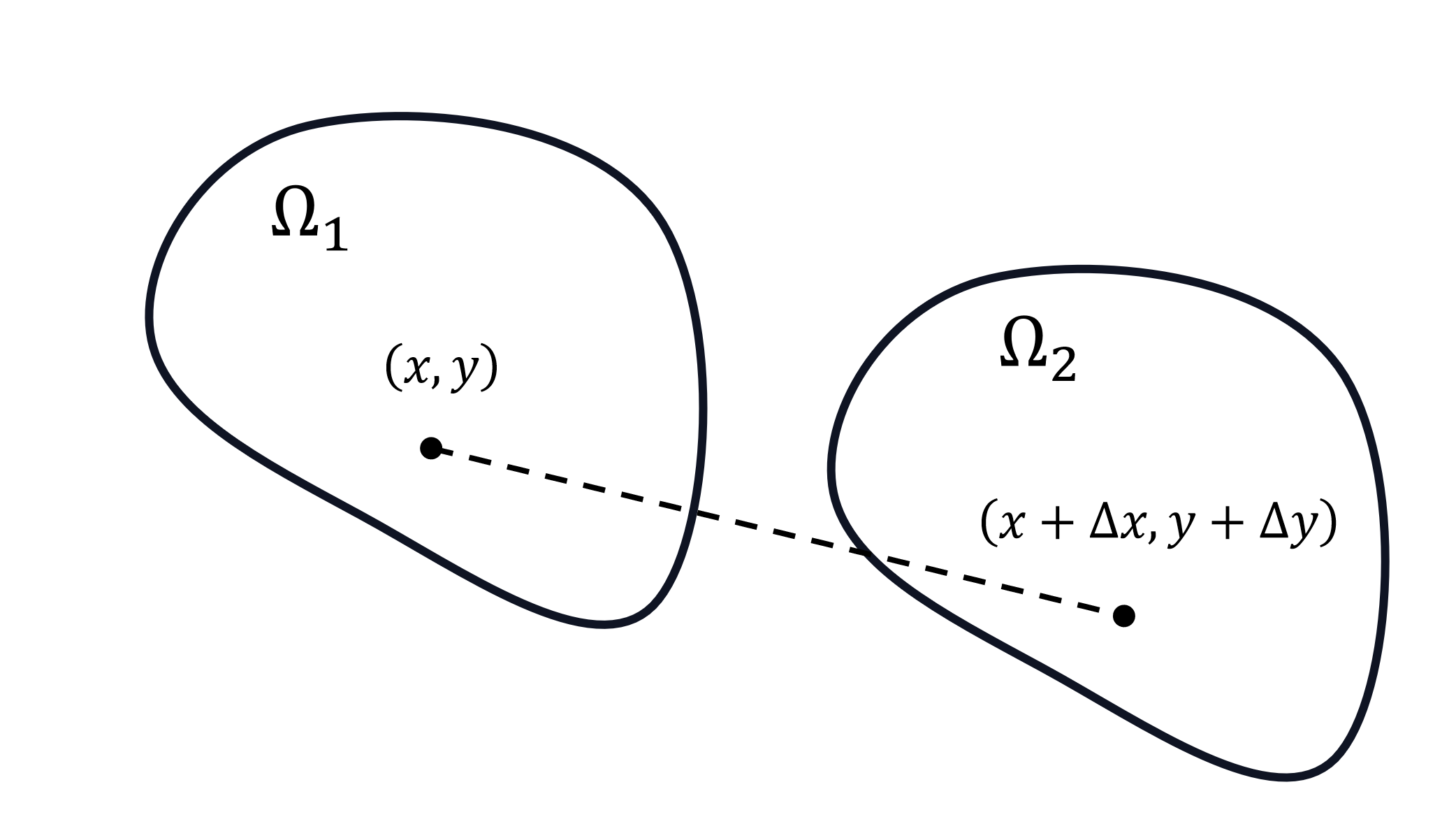}
\caption{Schematic diagram of moving $\Omega_1$ to $\Omega_2$.}
\label{Fig_DS}
\end{figure}

Next, we give a brief analysis of the effectiveness of the separation strategy.

Taking a 2D model as an example, assume that $u_1$ and $u_2$ are functions on the subdomains $\Omega_1$ and $\Omega_2$, respectively, satisfying Eqs. \eqref{2.1} and \eqref{bdc}, and that the diffusion coefficient $\kappa(x,y)$ , the source term $Q(x,y)$ and the boundary condition $g(x,y)$ are equal to the values of the corresponding position after moving the computational domain from $\Omega_1$ to $\Omega_2$ , as shown in Figure \ref{Fig_DS}, i.e.,
\begin{align} \label{cond13}
    \kappa(x,y) &=\kappa(x+\Delta x,y+\Delta y) , \quad (x,y) \in \Omega_1, (x+\Delta x,y+\Delta y)\in \Omega_2, \nonumber\\
    Q(x,y) &=Q(x+\Delta x,y+\Delta y) ,\quad (x,y) \in \Omega_1,
    (x+\Delta x,y+\Delta y)\in \Omega_2,\\ 
    g(x,y) &=g(x+\Delta x,y+\Delta y),\quad (x,y) \in \partial \Omega_1,(x+\Delta x,y+\Delta y) \in \partial \Omega_2.\nonumber
\end{align}
Then we have 
\begin{equation}
    u_1(x,y)= u_2(x+\Delta x,y+\Delta y),\quad (x,y) \in \Omega_1.
\end{equation}
This is an obvious result, and we can also give a simple proof.\\
Let 
\begin{equation}
    w(x,y) = u_1(x,y) - u_2(x+\Delta x,y+\Delta y),
\end{equation}
and then $w(x,y)$ satisfies
\begin{equation}\label{eq:w}
\begin{split}
\begin{cases}
     -\nabla\cdot\kappa(x,y)\nabla w(x,y) = 0, 
     \quad (x,y) \in \Omega_1,\\
     \left.w(x,y)\right|_{(x,y)\in\partial\Omega_1} = 0.
\end{cases}
\end{split}
\end{equation}
According the extremum principle, we get $w(x,y)\equiv0$. 

The above analysis shows that after translating the computational domain $\Omega_1$ to a new location $\Omega_2$, if the condition \eqref{cond13} is satisfied, the solutions of both domains have the same structure for Eqs.~\eqref{2.1}-~\eqref{bdc}.

\begin{remark}
Domain separation is a constructive strategy that makes it possible to express a class of non-smooth functions with a single neural network, fully exploiting the mesh-free advantage of the PINN method. 
Compared to the conventional method of using multiple neural networks to solve multi-material diffusion problems, the DS-PINN method using this strategy is not only easy to implement, but also does not require iteration between networks, resulting in relatively high computational efficiency.
\end{remark}

\subsection{Adding the special term representing the interface connection condition to the loss function.}

For multi-material diffusion problems, the solution at the interface is critical, and obtaining a highly accurate numerical solution near the interface is a very challenging task for any kind of numerical method.
Since $u$ has no second-order partial derivatives at the interface, the standard PINN method based on the automatic differentiation technique cannot incorporate the sampled points at the interface as residuals in the training of the neural network, which inevitably leads to inaccurate solutions near the interface. Sect.~\ref{linkcond} gives continuity conditions that the solution at the interface should satisfy. Introducing this connection condition into PINN can fill the gap of missing information at the interface. Furthermore, by introducing a domain separation strategy, Sect.~\ref{domsep} overcomes the problem of two derivative values at one location caused by $u$ not being continuously differentiable at the interface,  which is another hurdle
for the standard PINN method to solve the heterogeneous diffusion problem.

Based on the work in the previous two sections, we can easily improve the standard PINN by adding the interface continuity condition to the loss function, so that the interface information is introduced into the neural network.
In return, the prediction function we obtain can give very accurate predictions near the interface. For simplicity, we will only discuss the case of a single interface for two materials. However, the case of multiple interfaces for multiple materials can be treated in a similar manner.

By adding the \emph{interface loss term} in the loss function, the loss function of the DS-PINN is given by 
\begin{equation}\label{loss_DS}
\mathcal{L}(\theta;\Sigma) = {w_{b}}\mathcal{L}_{b}(\theta;\tau_{b})+
{w_{r}} \mathcal{L}_{r}(\theta;\tau_{r})
+{w_{\Gamma}} \mathcal{L}_{\Gamma}(\theta;\tau_{\Gamma}),
\end{equation}
where 
\begin{align}
\mathcal{L}_{b}(\theta;\tau_{b}) =& \frac{1}{N_b}\sum_{i=1}^{N_b}\left|u_\theta(X'_i)-g(X_i)\right|^2,\label{DS_PINNLb}
\end{align}
\begin{align}
\mathcal{L}_{r}(\theta;\tau_{r})=&\frac{1}{N_r}\sum_{i=1}^{N_{r}}\left|-\nabla\cdot \kappa(X_i)\nabla u(X'_i) - Q\left(X_i\right)
\right|^2,  \label{DS_PINNLr}
\end{align}
\begin{align}
\mathcal{L}_\Gamma(\theta;\tau_{\Gamma})
=& \frac{1}{N_\Gamma}\sum_{i=1}^{N_\Gamma}\left(\left|u_\theta(X_i^{(1)})-u_\theta(X_i^{(2)})\right|^2+\right.\notag\\
&~\left.
\left|
-\kappa_1(X_i^{(1)})\nabla u(X_i^{(1)})\cdot\boldsymbol{n}_1-\kappa_2(X_i^{(2)})\nabla u(X_i^{(2)})\cdot\boldsymbol{n}_2
\right|^2\right).\label{loss_gamma}
\end{align}
The training data set
$\Sigma\left(\tau_{b},\tau_{r}, \tau_{\Gamma}\right)$ is as follows:
\begin{align}
&\tau_b=\{\left(X_i,g(X_i)\right)|X_i\in\partial\Omega\}_{i=1}^{N_b},\\
&\tau_{r}=\{X_i \in\Omega\}_{i=1}^{N_{r}},\\
&\tau_{\Gamma}=\{X_i|X_i\in \Gamma\}
_{i=1}^{N_\Gamma}.\label{GMPT}
\end{align}

In Eqs. \eqref{DS_PINNLb} and \eqref{DS_PINNLr}, $X_i$ represents a sampling point on the boundary of the original domain $\partial\Omega$ or in the original domain $\Omega$, and $X'_i$ is the matching point of $X_i$. 
If the subdomain containing $X_i$ has not moved, then $X'_i=X_i$;
if the subdomain containing $X_i$ has moved by a distance of $\Delta L$, then $X'_i=X_i+\Delta L$.

In Eqs. \eqref{loss_gamma} and \eqref{GMPT}, $X_i$ represents a sampling point on the material interface $\Gamma$ in the original domain $\Omega$, and $\Gamma$ are referred to as $\Gamma_1$ and $\Gamma_2$ in two neighboring subdomains $\Omega_1$  and $\Omega_2$, respectively.  
$X^{(1)}_i$ and $X^{(2)}_i$ are two matching points of  $X_i$ belonging to $\Gamma_1$ and $\Gamma_2$. 
If the subdomain containing $X^{(k)}_i$ has not moved, then $X^{(k)}_i=X_i$, k=1 or 2;
if the subdomain containing $X^{(k)}_i$ has moved by a distance of $\Delta L$, then $X^{(k)}_i=X_i+\Delta L$,  k=1 or 2.

Note that, according to  Eq. \eqref{cond13}, we use $g(X_i)$, $\kappa(X_i)$ and $Q(X_i)$ to replace with $g(X'_i)$,$\kappa(X'_i)$ and $Q(X'_i)$ in Eqs. \eqref{DS_PINNLb} and \eqref{DS_PINNLr}.
In addition, if higher accuracy is needed at the material interface, one can increase $\omega_\Gamma$.

So far, we have obtained the loss function for the multi-material diffusion equation, and by optimizing this loss function, we can get the prediction function. The last task is to map the prediction function of the subdomain shifted by a certain distance back to its original position using the following formula
\begin{equation} \label{MapU}
   u_\theta(X)=u_\theta(X'), \quad X \in \bar{\Omega}.
\end{equation}
Similarly, $X'$ is the matching point of $X$.
Depending on whether the subdomain containing $X$  is moved, $X'= X$ or $X'= X+\Delta L$.

\subsection{Normalizing loss terms to improve the training performance}

The loss function of the standard PINN method contains two main terms, the supervised loss term and the residual loss term, whose variables usually have different physical meanings and are of different magnitudes, so that combining them for optimization will generally result in the numerically smaller term not being reasonably optimized, 
leading to the final prediction deviating from the reference solution. 
The question of how to balance the different terms in the loss function plays a key role in the PINN method, and some researchers have made important progress \cite{HE1,Xie,Wu,Wight}. 

It should be noted that normalizing each loss term in the loss function according to the characteristics of the equation not only facilitates the implementation of the optimization algorithm, but also helps to balance the importance of each loss term, eliminates poor training results due to different orders of magnitude, and improves the computational accuracy of the prediction function.
Based on the governing equation \eqref{2.1} and its boundary condition \eqref{bdc}, a strategy for normalizing the supervised term \eqref{DS_PINNLb} and the residual term \eqref{DS_PINNLr} is given below.

Considering that $g(X)$ and $Q(X)$ may reflect the magnitude of $u$ and $u_{xx}+u_{yy}$, respectively, let the normalization factors $\zeta_b$ and $\zeta_r$ are as follows
\begin{align}
\zeta_b&=\frac{1}{N_b}\sum_{i=1}^{N_b}\left|g(X_i)\right|^2,
\label{DS_PINNLb_factor}\\
\zeta_r&=\frac{1}{N_r}\sum_{i=1}^{N_{r}}\left| Q\left(X_i\right)
\right|^2.  \label{DS_PINNLr_factor}
\end{align}
Next, we define the new loss terms
\begin{equation}\label{normalDS_PINNLb}
\mathcal{\Tilde{L}}_{b}(\theta;\tau_{b}) = 
\begin{cases}
\frac{1}{\zeta_b}\mathcal{L}_{b}(\theta;\tau_{b}),&~\text{if}~ \zeta_b \neq 0,\\
\mathcal{L}_{b}(\theta;\tau_{b}),&~\text{if}~ \zeta_b = 0.
\end{cases}
\end{equation}

\begin{equation}\label{normalDS_PINNLr}
 \mathcal{\Tilde{L}}_{r}(\theta;\tau_{r}) = 
\begin{cases}
\frac{1}{\zeta_r}\mathcal{L}_{r}(\theta;\tau_{r}),&~\text{if}~ \zeta_r \neq 0,\\
\mathcal{L}_{r}(\theta;\tau_{r}),&~\text{if}~ \zeta_r = 0.
\end{cases}
\end{equation}
Note that, although the normalization factors cannot normalize the value of each loss term to $[0,1]$, it is able to largely eliminate the effect of the magnitude. Now, the loss function can be rewritten as follow
\begin{equation}\label{normlizeloss_DS}
\mathcal{L}(\theta;\Sigma) = {w_{b}}\mathcal{\Tilde{L}}_{b}(\theta;\tau_{b})+
{w_{r}} \mathcal{\Tilde{L}}_{r}(\theta;\tau_{r})
+{w_{\Gamma}} \mathcal{L}_{\Gamma}(\theta;\tau_{\Gamma}),
\end{equation}
and we call the PINN method using this loss function the \emph{normalized DS-PINN}, denoted by \emph{nDS-PINN}.

\begin{remark}
The interface loss term $\mathcal{L}_{\Gamma}(\theta;\tau_{\Gamma})$ isn't normalized because it reflects the continuity conditions ~\eqref{cu_new}-\eqref{cFlux_new}.
This is similar to the case where $\zeta_b=0$ or $\zeta_r=0$ in Eq. \eqref{normalDS_PINNLb} or Eq. \eqref{normalDS_PINNLr}.
The interface connection condition can be generalized to the general form
\begin{align}
[\![u(X)]\!]_{\Gamma}&=\Phi(X),\label{cu_gen}\\
[\![-\kappa(X)\nabla u(X)\cdot\boldsymbol{n} ]\!]_{\Gamma}&=\Psi(X).\label{cFlux_gen}
\end{align}
The continuity conditions \eqref{cu_new} and \eqref{cFlux_new} are considered as a special case where $\Phi(X)=\Psi(X)=0$.
For a class of interface problems with jump conditions, i.e., for the case where
$\Phi(X)\neq0$ and $\Psi(X)\neq0$, 
we can use $\Phi(X)$ and $\Psi(X)$ 
to normalize 
$\mathcal{L}_{\Gamma}(\theta;\tau_{\Gamma})$
in the same way as Eqs. \eqref{DS_PINNLb_factor}-\eqref{normalDS_PINNLr}.
\end{remark}

\section{Numerical experiments}\label{NEX}
In this section, we give several numerical experiments to demonstrate the performance of the proposed method in this work. 
In Sect.~\ref{Ex1}, we test the performance of the new methods by solving a typical two-material diffusion equation, and we show the results based on different separation distances.
In Sect.~\ref{Ex3}, we solve a multi-material diffusion equation with 4 subdomains. 
In Sect.~\ref{Ex4}, we present a diffusion model with a special computational domain. The computational domain contains different materials within a circular subdomain at its center.
In Sect.~\ref{Ex5}, we test the ability of the new method for diffusion problems with jump conditions at the interface.



We use the deep learning framework TensorFlow(version 1.5) to develop the code. The datatype of all variables is  {\itshape float32}. For all the numerical experiments, we use the Adam optimizer to run 2000 iterations and then switch to the L-BFGS optimizer until convergence. All parameters and termination criteria of the L-BFGS optimizer are considered as proposed   in Ref.~\cite{Byrd}. 
Before training, the parameters of the neural network are randomly initialized by using the Xavier scheme~\cite{XuC}.
All numerical experiments use the deep neural network with 5 hidden layers of 50 neurons each.  

The accuracy of the trained model is evaluated using the relative $\mathbb{L}_2$ error, which is defined as follows:
\begin{equation}
\left\|e\right\|_{\mathbb{L}_2}=\frac{\sqrt{\sum_{i=1}^N\left|u_\theta(X_i)-u^*(X_i)\right|^2}}{\sqrt{\sum_{i=1}^N\left|u^*(X_i)\right|^2}},
\end{equation}
where we suppose
$\sum_{i=1}^N\left|u^*(X_i)\right|^2 \neq 0$,
$u^*(X_i)$ is the exact solution or the reference solution, and $u_\theta(X_i)$ is the neural network prediction for $N=10000$ test points which are uniformly distributed over the computational domain. 
For the number of training points in this section, we set $N_f=5000$ for each subdomain, $N_b=500$ for each boundary, and $N_{\Gamma}=2000$ for each interface.

\begin{remark}
Different hyperparameters, such as the model architecture, the size of the training data, the weights, the optimizer, etc., can cause the PINN method to produce different computational results. To test the robustness and effectiveness of our method, we use the same set of parameters for all numerical examples in this paper. We believe that finer tuning of the parameters could yield even better results.
\end{remark}

\subsection{Typical two-material diffusion problems}\label{Ex1}

Consider the following diffusion problem with two materials in the computational domain:
\begin{equation}\label{StandardEq}
    \begin{cases}
        -\nabla\cdot\left(\kappa(x,y)\nabla u\right)=Q(x,y),&~(x,y)\in\Omega=(0,1)\times(0,1),\\
    u(x,y)=0,&~(x,y)\in\partial\Omega,
    \end{cases}
\end{equation}
where
\begin{equation}
    \kappa(x,y)=
    \begin{cases}
        4,~(x,y)\in\left(0,\frac{2}{3}\right]\times(0,1),\\
        1,~(x,y)\in\left(\frac{2}{3},1\right)\times(0,1),
    \end{cases}
\end{equation}
and
\begin{equation}
    Q(x,y)=
    \begin{cases}
        20\pi^2\sin\pi x\sin 2\pi y,~(x,y)\in\left(0,\frac{2}{3}\right]\times(0,1),\\
        20\pi^2\sin4\pi x\sin 2\pi y,~(x,y)\in\left(\frac{2}{3},1\right)\times(0,1).
    \end{cases}
\end{equation}
The exact solution of this equation is
\begin{equation}
    u(x,y)=
    \begin{cases}
        \sin\pi x\sin 2\pi y,~(x,y)\in\left(0,\frac{2}{3}\right]\times(0,1),\\
        \sin4\pi x\sin 2\pi y,~(x,y)\in\left(\frac{2}{3},1\right]\times(0,1).
    \end{cases}
\end{equation}

\begin{figure}[h]\label{fig:case1_exact}
\centering
\includegraphics[width=0.65\textwidth]{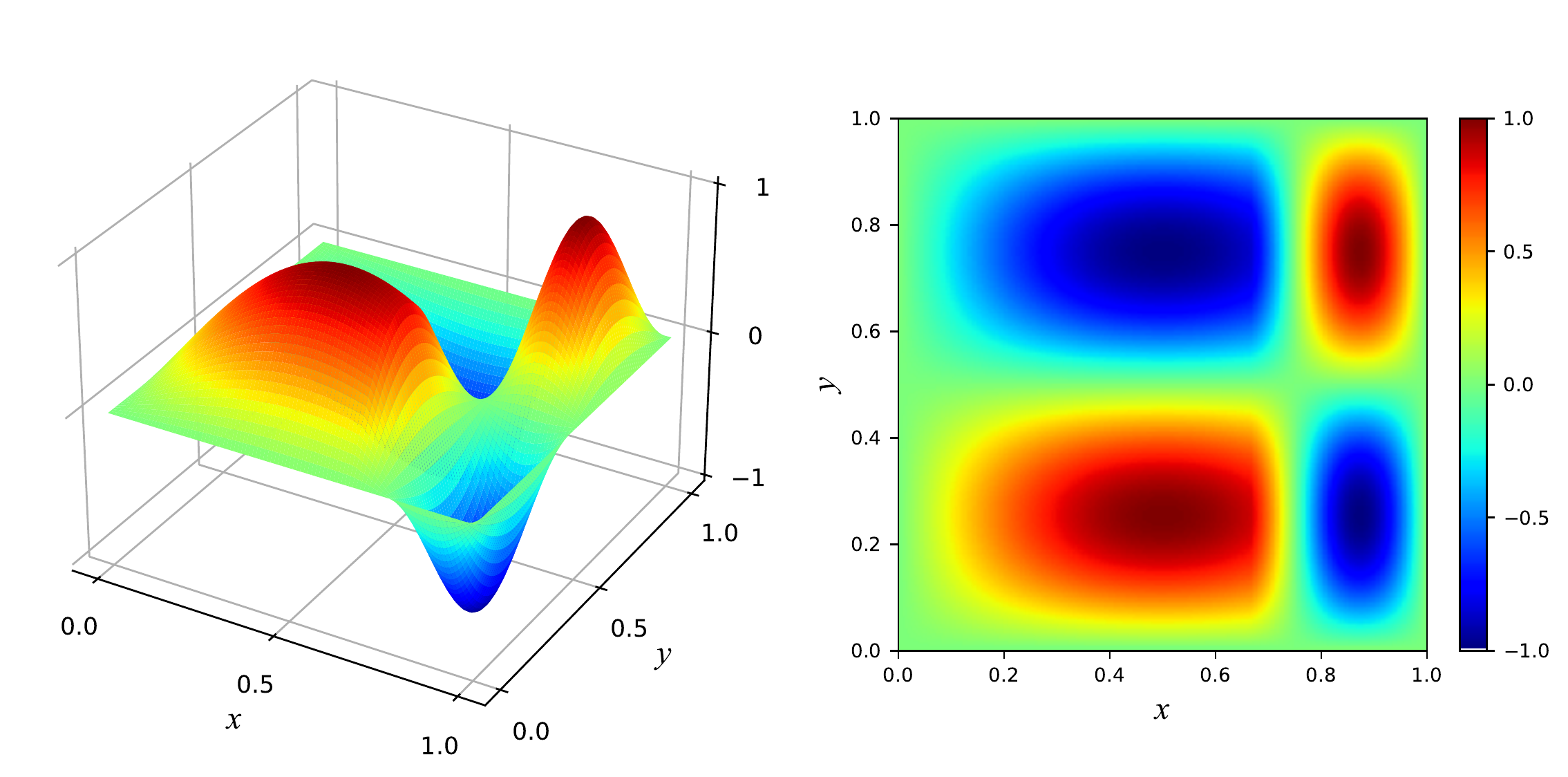}
\caption{ The exact solution of Sect.~\ref{Ex1}. }
\end{figure}

Figure \ref{fig:case1_exact} shows the exact solution of this example, expressed in z-coordinate value and color, respectively.
The solution $u(x,t)$ is continuous at the interface $x=\frac{2}{3}$, but its partial derivative $u_x$ is discontinuous at the interface. 
\begin{figure}[h]
\centering
\includegraphics[width=0.9\textwidth]{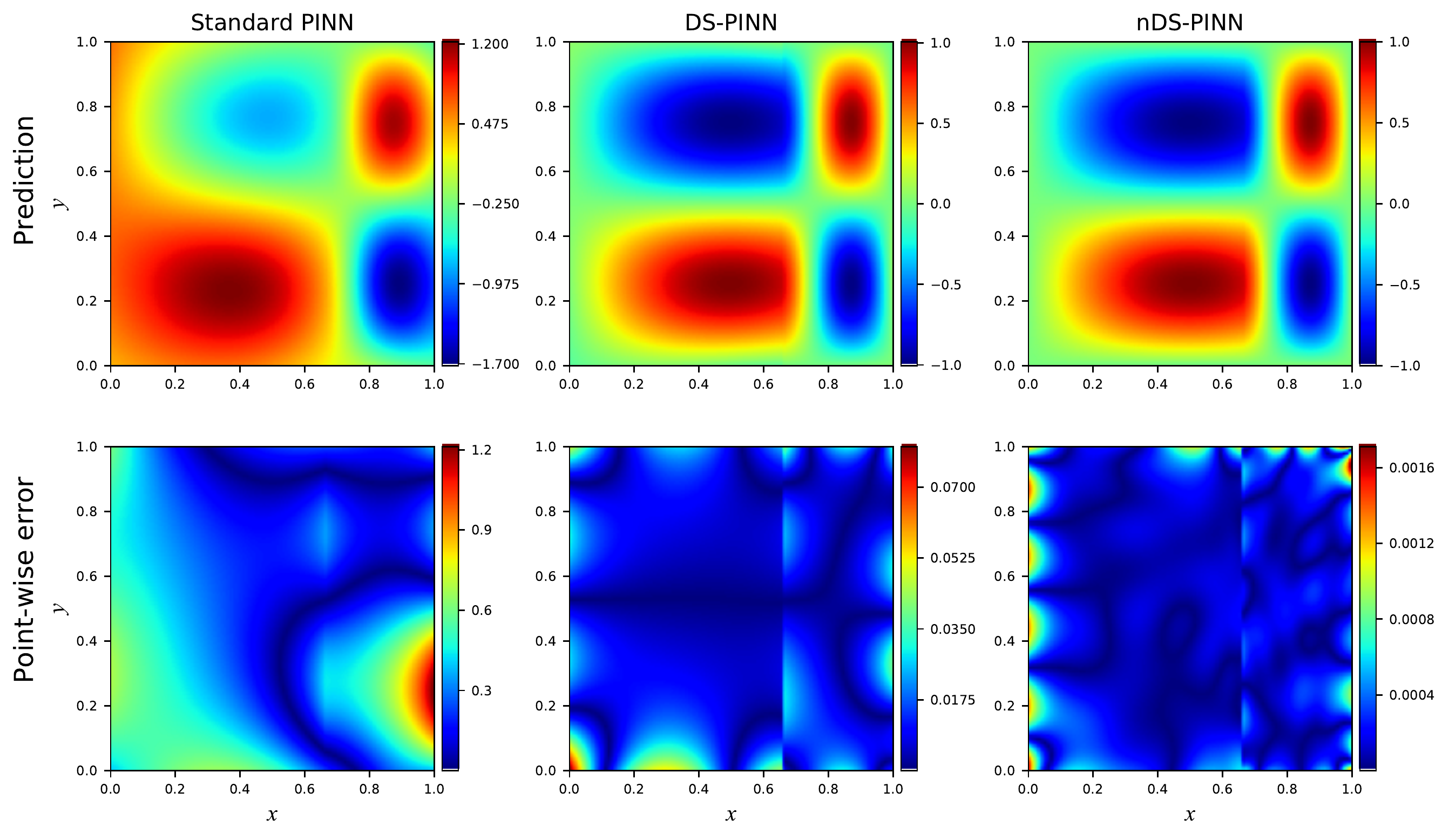}
\caption{Predictions and point-wise errors of three PINN methods for Example \ref{Ex1} }\label{case-twomat}
\end{figure}

For this model, we set the domain separation distance ${d}=0.1$. 
Table \ref{tabeforex1} shows the relative $\mathbb{L}_2$ error of different methods.
Figure \ref{case-twomat} shows the the prediction and the point-wise error of different PINN methods.
It can be seen that the standard PINN method gives a smooth prediction, and due to the lack of interface information, its prediction is of low accuracy.
However, the DS-PINN method provides a highly accurate prediction, and by normalizing the residual terms, the nDS-PINN method further improves the prediction accuracy by more than an order of magnitude.
\begin{table}[h]\label{tabeforex1}
		\setlength{\abovecaptionskip}{0cm}
		\setlength{\belowcaptionskip}{0.2cm}
	\centering
	\caption{$\left\|e\right\|_{\mathbb{L}_2}$ errors of three PINN methods for Sect.~\ref{Ex1}.}
	\begin{tabular}{c|c}
		\hline
		Method  & $\left\|e\right\|_{\mathbb{L}_2}$ \\ 
		 \hline
		 Standard PINN        &   $6.0\pm0.4\times10^{-1}$\\
		 DS-PINN (this work)      &    $4.1\pm1.3\times10^{-2}$\\
		 nDS-PINN (this work)    &   $7.3\pm1.6\times10^{-4}$\\
		 \hline
	\end{tabular}
\end{table}

\begin{figure}[h]
\centering
\includegraphics[width=0.6\textwidth]{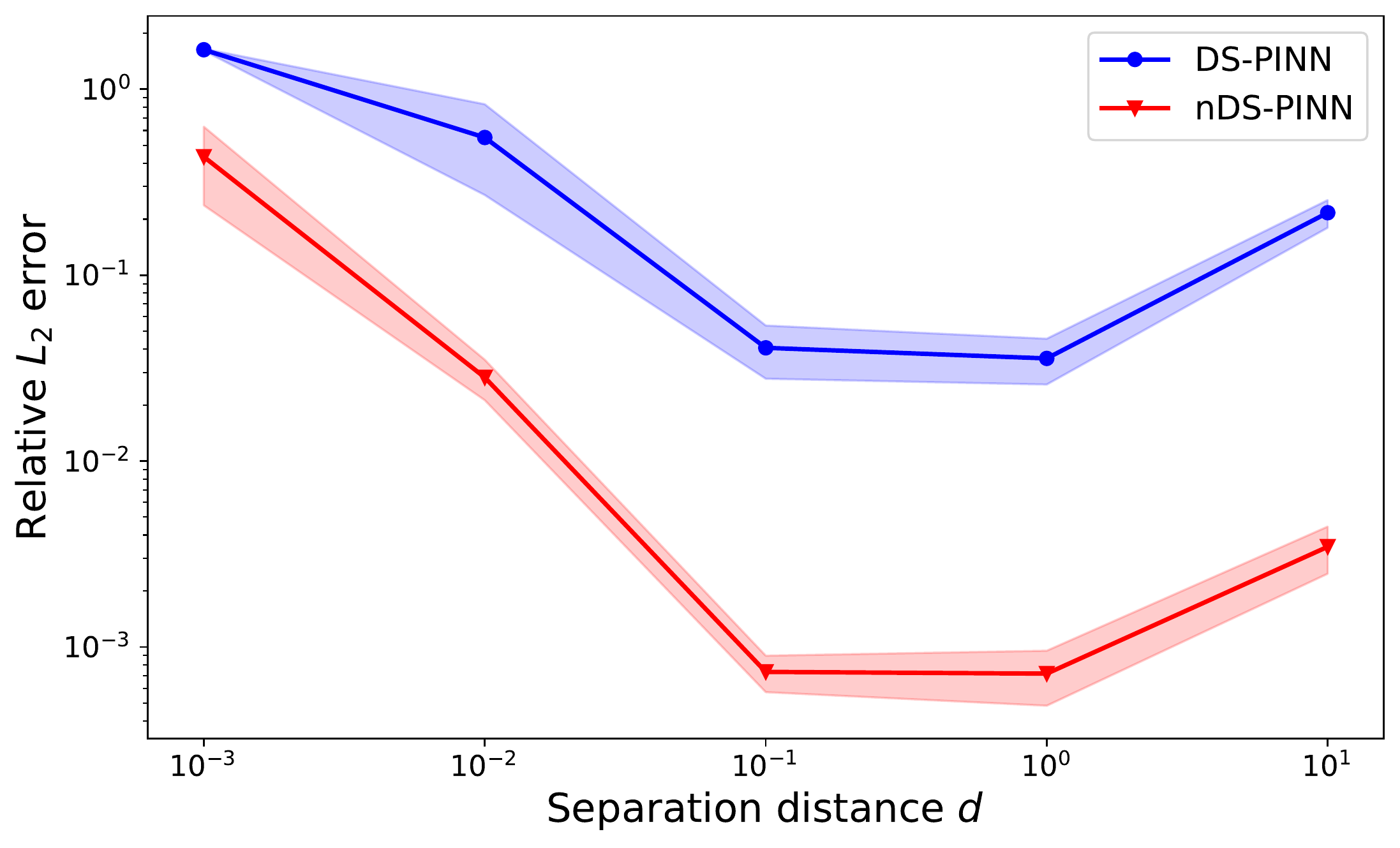}
\caption{$\left\|e\right\|_{\mathbb{L}_2}$ errors of the two new methods with different separation distances(obtained from 5 independently repeated experiments).
}\label{fig:case1_delta}
\end{figure}

How does the separation distance between subdomains affect the solution accuracy?
To investigate this question, we use different separation distances to test this computational model, and the figure \ref{fig:case1_delta} 
shows the effect of different ${d}$ on the prediction accuracy.

As can be seen from the Figure \ref{fig:case1_delta}, if ${d}$ is too small (${d}<10^{-2}$), the error is very large.
The reason is that, for a small ${d}$, when the diffusion coefficients on both sides of the material interface differ greatly, it is equivalent to solving a problem with a large variation of the derivative over a narrow interval, which is computationally very difficult and therefore yields inaccurate results without special tricks. 
Conversely, if ${d}$ is too large (${d}>10$), 
then the range expressed by the neural network increases significantly, 
and moreover, the gaps between the subdomains do not assign the sampling points and do not participate in the training of the neural network, and so they belong to the invalid computational domain. Obviously, if the proportion of invalid regions is too large, the accuracy of the prediction function will inevitably be low.

\begin{remark}
Regarding the standard PINN method, although one can get a higher accuracy by changing some hyper-parameters, its prediction accuracy near the material interface does not change significantly. 
This remark is also adapted to the discussion of the standard PINN method in the later examples.
\end{remark}

\begin{remark}
Choosing an appropriate separation distance is beneficial for improving the computational accuracy of the model.
We believe that ${d}$ should be chosen based on the fact that it should be proportional to $\left|\kappa_1-\kappa_2\right|$, while keeping the percentage of added invalid computational area as small as possible. Here $\kappa_1$ and $\kappa_2$ represent the diffusion coefficients on both sides of the interface. 
\end{remark}

\subsection{Multi-material diffusion problems}\label{Ex3}
In this subsection, we examine a multi-material diffusion example and the separation distance $d=0.1$.
The governing equation is the same as Eq.~\eqref{StandardEq}.
Suppose that the computational domain $\Omega$ consists of 4 subdomains with different diffusion coefficients,
\begin{equation}\label{MultMkapp}
    \kappa(x,y)=
    \begin{cases}
        4,&~(x,y)\in\left(-1,0\right]\times(-1,0],\\
        1,&~(x,y)\in\left(0,1\right)\times(-1,0],\\
        2,&~(x,y)\in\left(0,1\right)\times(0,1),\\
        1,&~(x,y)\in\left(-1,0\right]\times(0,1).
    \end{cases}
\end{equation}
We also assume that this problem has the exact solution as follows:
\begin{equation}\label{multiexact}
    u(x,y)=
    \begin{cases}
        \sin\pi x\sin\pi y,&~(x,y)\in\left[-1,0\right]\times[-1,0],\\
        4\sin\pi x\sin\pi y,&~(x,y)\in\left(0,1\right]\times[-1,0],\\
        2\sin\pi x\sin\pi y,&~(x,y)\in\left(0,1\right]\times(0,1],\\
        4\sin\pi x\sin\pi y,&~(x,y)\in\left[-1,0\right]\times(0,1].
    \end{cases}
\end{equation}
To test the performance of the new methods, we solve this model using DS-PINN and nDS-PINN. The source term $Q(x,y)$ and the boundary conditions are derived from the exact solution \eqref{multiexact}.

\begin{figure}[!htbp]
\centering
\includegraphics[width=0.5\textwidth]{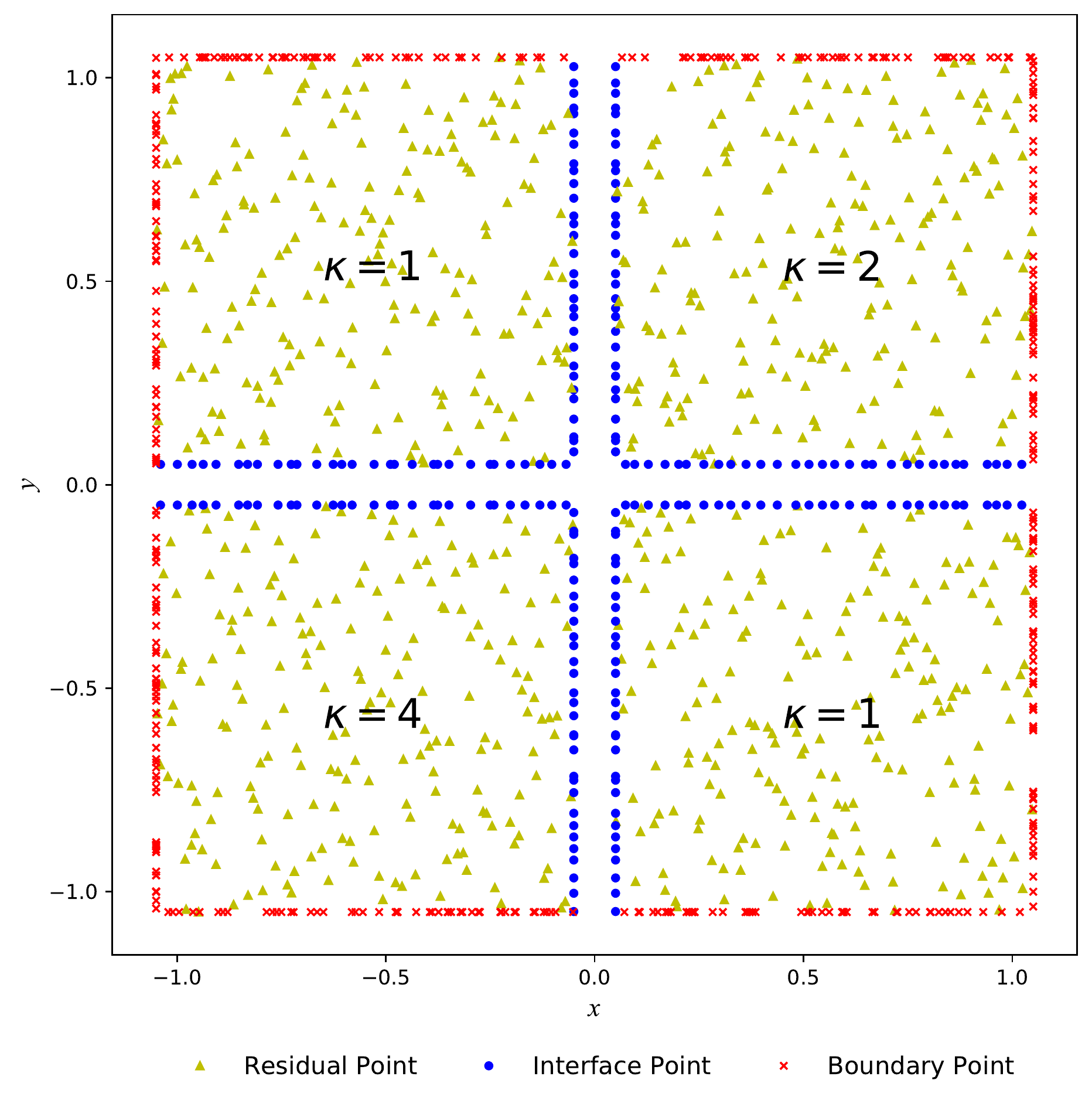}
\caption{Schematic diagram of domain separation and training point layout 
of Sect.~\ref{Ex3}.}
\label{fig:case3_points}
\end{figure}
Figure \ref{fig:case3_points} shows a schematic diagram of the sampling of training points in different subdomains, and the training points consist of three types, including residual points located inside the subdomains, supervised points located at the boundaries, and interface points located at the material interfaces.

This example is somewhat complicated. If multiple neural networks were conventionally used to handle the material interface, this problem would require four neural networks, which would not only be difficult to implement, but also computationally inefficient. With our methods, this problem can be easily solved with only a single neural network.
\begin{figure}[!htbp]
\centering
\includegraphics[width=0.9\textwidth]{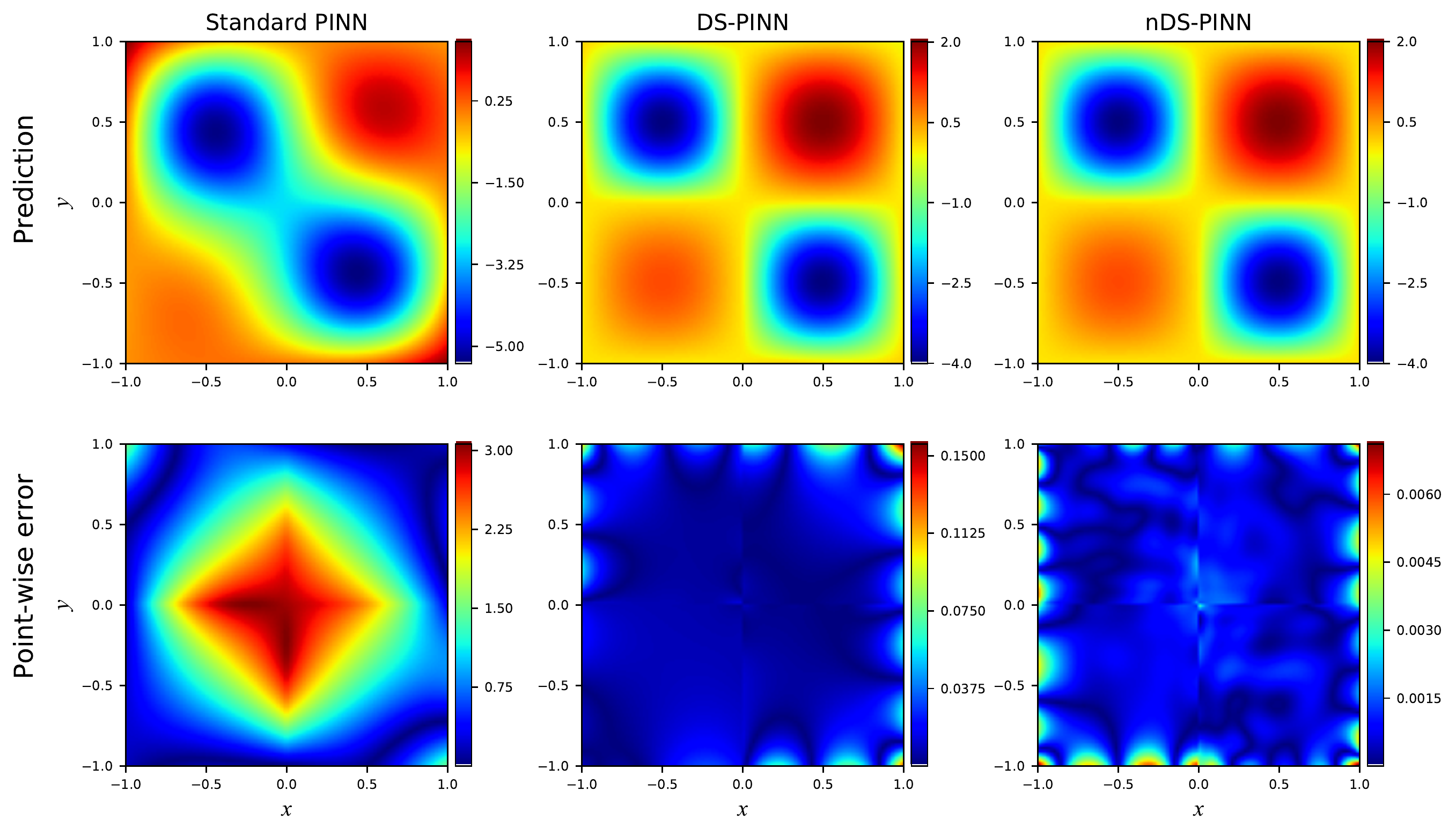}
\caption{Predictions and point-wise errors of three PINN methods for Sect.~\ref{Ex3} }\label{case-multimat}
\label{fig:multmatexample}
\end{figure}
\begin{table}[!htbp]\label{tabeforex2}
		\setlength{\abovecaptionskip}{0cm}
		\setlength{\belowcaptionskip}{0.2cm}
	\centering
	\caption{$\left\|e\right\|_{\mathbb{L}_2}$ errors of three PINN methods for Sect.~\ref{Ex3}.}
	\begin{tabular}{c|c}
		\hline
		Method  & $\left\|e\right\|_{\mathbb{L}_2}$ \\ 
		 \hline
		 Standard PINN        &   $8.96\pm0.03\times10^{-1}$\\
		 DS-PINN (this work)      &    $1.5\pm0.2\times10^{-3}$\\
		 nDS-PINN (this work)    &   $5.2\pm1.6\times10^{-4}$\\
		 \hline
	\end{tabular}
\end{table}

Table \ref{tabeforex2} shows the relative $\mathbb{L}_2$ error of different methods, and
Figure \ref{fig:multmatexample} shows the the prediction and the point-wise error of different PINN methods.
Similar to the results of the previous example, the standard PINN gives a poor prediction for this model, while the DS-PINN gives a satisfactory prediction and the nDS-PINN gives an accurate prediction.

\subsection{The diffusion problem with the heterogeneous material located inside the computational domain}\label{Ex4}

In practical applications such as heat transfer and oil reservoir simulation, it is common for a material to be completely enveloped by another material. The purpose of this section is to test the ability of our method to handle this case.

Consider the problem as follows:
\begin{equation}
\begin{cases}
     -\nabla\cdot\left(\kappa(x,y)\nabla u\right)=Q(x,y),&~(x,y)\in\Omega = (-2,2)\times(-2,2),\\
    u(x,\pm 2)=\sin(\frac{\pi}{4}(x^2+3)),&~x\in [-2,2],\\
    u(\pm 2,y)=\sin(\frac{\pi}{4}(y^2+3)),&~y\in [-2,2],   
\end{cases}
\end{equation}
where
\begin{equation*}
    \kappa(x,y)=
    \begin{cases}
        1,~ & (x,y)\in\Omega_1 =\{(x,y)|x^2+y^2<1\},\\
        4,~ &(x,y)\in \Omega \backslash \Omega_1,
    \end{cases}
\end{equation*}
and the source term 
\begin{equation*}
    Q(x,y)=
    \begin{cases}
        -4\pi\cos\left(\pi(x^2+y^2-1)\right)+4\pi^2(x^2+y^2)\sin\left(\pi(x^2+y^2-1)\right),\\ \hfill(x,y)\in \Omega_1,\\
        -4\pi\cos\left(\frac{\pi}{4}(x^2+y^2-1)\right)+\pi^2(x^2+y^2)\sin\left(\frac{\pi}{4}(x^2+y^2-1)\right), \\ \hfill(x,y)\in\Omega\backslash \Omega_1.
    \end{cases}
\end{equation*}
The exact solution of this problem is
\begin{equation}
    u(x,y)=
    \begin{cases}
        \sin(\pi(x^2+y^2-1)),~ & (x,y)\in\Omega_1,\\
        \sin(\frac{\pi}{4}(x^2+y^2-1)),~ &(x,y)\in\Omega\backslash \Omega_1.
    \end{cases}
\end{equation}

\begin{figure}[!htbp]
\centering
\includegraphics[width=1.0\textwidth]{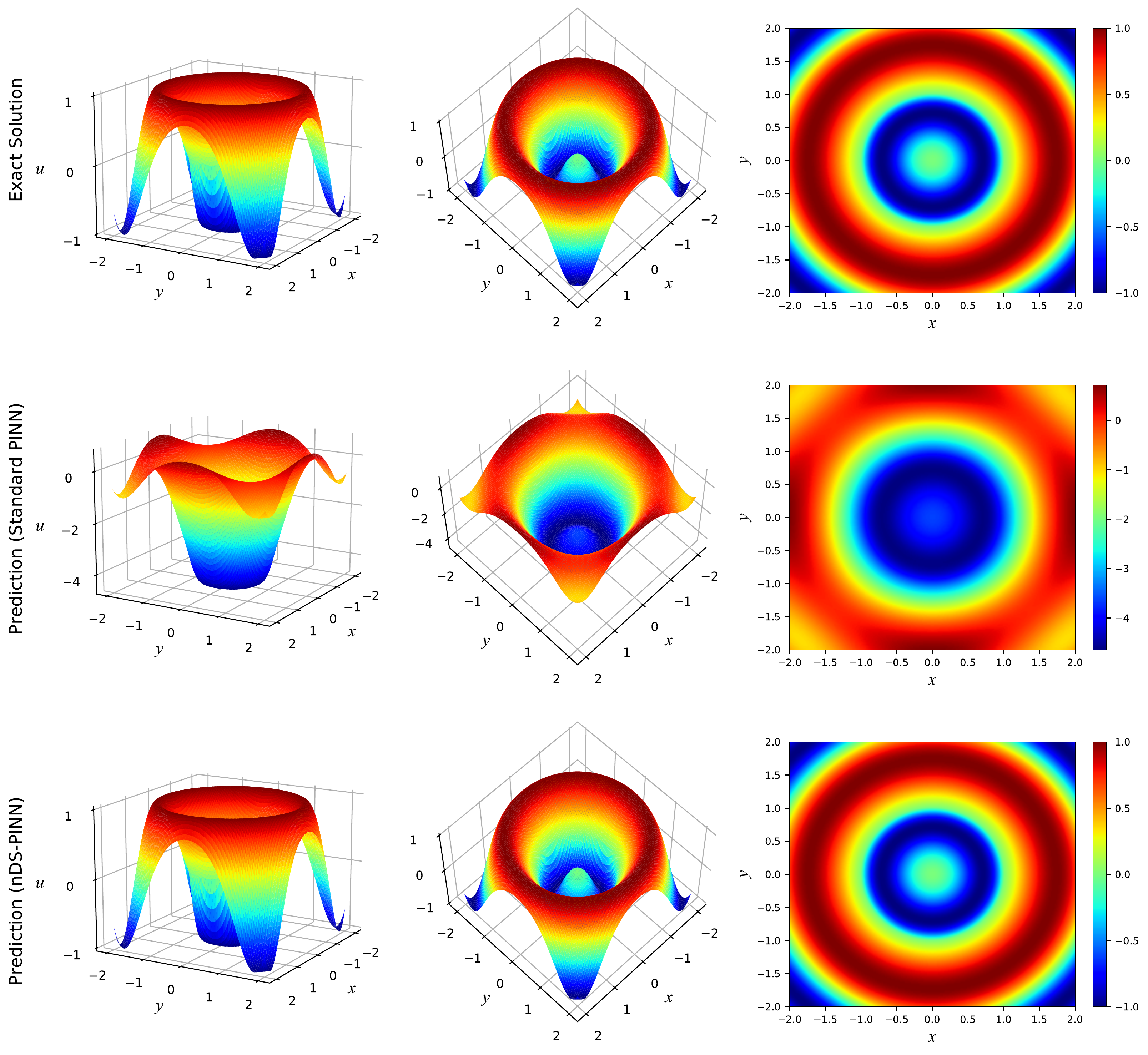}
\caption{ The exact solution of Sect.~\ref{Ex4} and the two predictions from the standard PINN and nDS-PINN.
\label{fig:bigpicforhole}
}
\end{figure}

\begin{figure}[!htbp]
\centering
\includegraphics[width=1.0\textwidth]{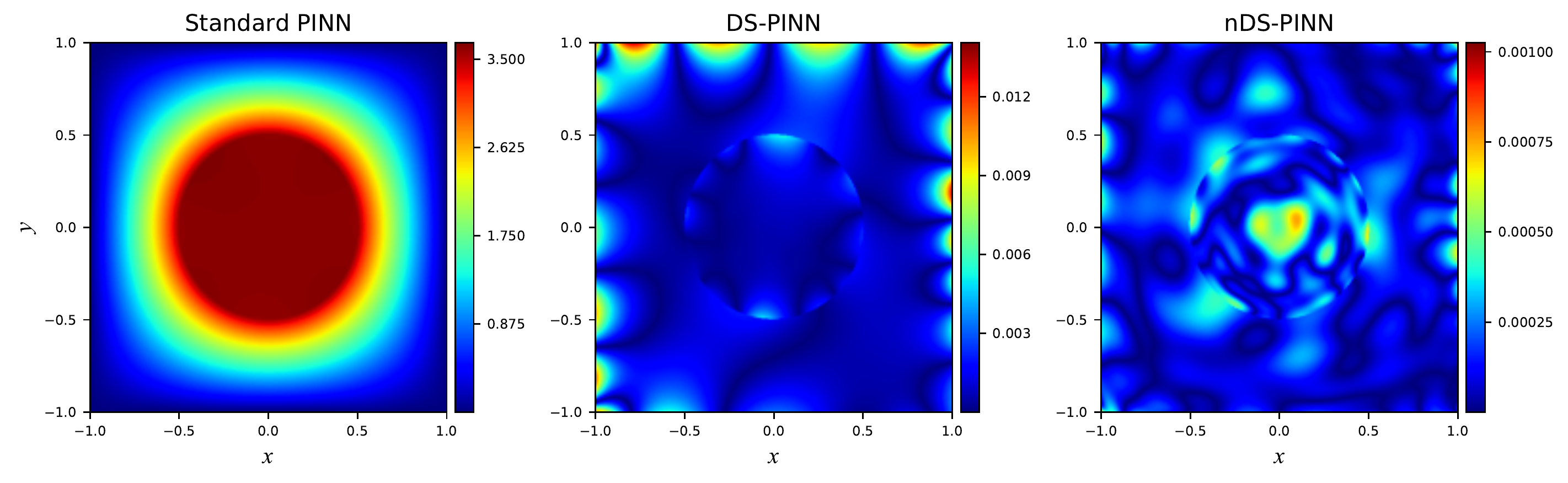}
\caption{Point-wise errors of three PINN methods for Sect.~\ref{Ex4}. 
\label{fig:bigpicforholepwerror}
}
\end{figure}

In this case, the separation distance $d=3.5$. Table \ref{examplehole} shows the results of three PINN method, and we can see that DS-PINN method and nDS-PINN method achieves satisfactory accuracy. The exact solution and the prediction using two PINN methods (the standard PINN method and the nDS-PINN method) are shown in Figure \ref{fig:bigpicforhole}, and the point-wise errors are shown in Figure \ref{fig:bigpicforholepwerror}.



\begin{table}[!htbp]\label{examplehole}
		\setlength{\abovecaptionskip}{0cm}
		\setlength{\belowcaptionskip}{0.2cm}
	\centering
	\caption{$\left\|e\right\|_{\mathbb{L}_2}$ errors of three PINN methods for Sect.~\ref{Ex4}.}
	\begin{tabular}{c|c}
		\hline
		Method  & $\left\|e\right\|_{\mathbb{L}_2}$ \\ 
		 \hline
		 Standard PINN        &   $2.96\pm0.01\times10^{0}$\\
		 DS-PINN (this work)      &    $3.3\pm1.3\times10^{-3}$\\
		 nDS-PINN (this work)    &   $3.3\pm1.4\times10^{-4}$\\
		 \hline
	\end{tabular}
\end{table}

It can be seen that the results of the standard PINN method are far from the exact solution. The main reason is that, due to the inability to arrange training points at the material interface, the exchange of information between $\Omega_1$ and 
$\Omega\backslash \Omega_1$ is hindered. 
At the same time, due to the fact that $\Omega_1$, which contains the heterogeneous material, is completely inside $\Omega$, its training is completely free from the constraints of the boundary conditions (the supervised term), which leads to a complete loss of control over the prediction of the central subdomain $\Omega_1$, and further leads to a complete deviation from the exact solution of the entire prediction.

On the other hand, The DS-PINN and nDS-PINN methods provide predictions that are highly consistent with the exact solution.
In particular, the nDS-PINN method with the normalization strategy shows a very high computational accuracy, both in terms of the relative $\mathbb{L}_2$ error and the point-wise error, which is one order of magnitude smaller than that of the DS-PINN method, showing excellent performance. Unlike the previous two examples, the boundary condition of this example is not zero, so not only the residual term is normalized, but also the supervised term is normalized.

\begin{figure}[!htbp]
\centering
\includegraphics[width=0.9\textwidth]{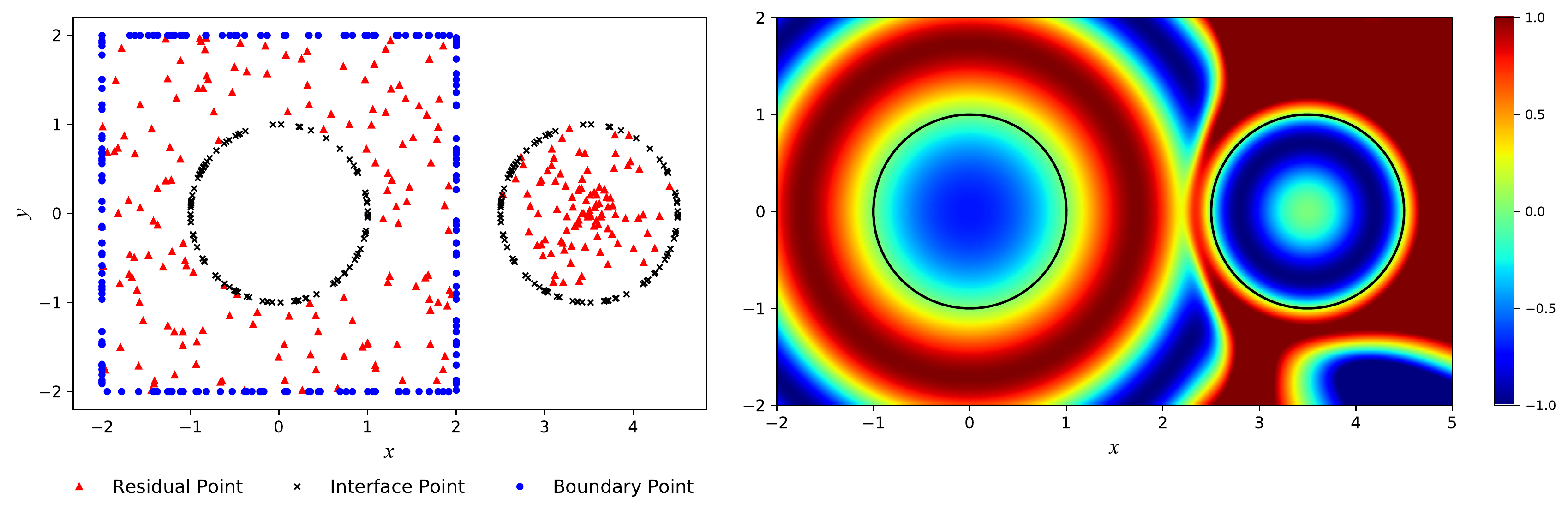}
\caption{Left: Schematic diagram of domain separation (${d}=3.5$) and training point layout
of Sect.~\ref{Ex4};
Right: The prediction on the extended  computational domain.}
\label{fig:case4_results_all}
\end{figure}
In Figure \ref{fig:case4_results_all}, the left image shows the schematic diagram of domain separation and it also gives the distribution of training points;
the right image shows the DS-PINN prediction for the whole extended domain after implementing the domain separation strategy. This is a very interesting picture, and we can see that the prediction on $\Omega_1$ and $\Omega \backslash \Omega_1$ match well with the exact solution at the corresponding locations.


\subsection{The diffusion problems with jump conditions at the interface}\label{Ex5}

In this paper we are mainly concerned with a class of multi-material diffusion problems formulated by Eqs.~\eqref{2.1}-\eqref{mmk}, for which the continuity conditions \eqref{cu_new} and \eqref{cFlux_new} should be satisfied at the material interface.  

However, there is a special class of heterogeneous diffusion problems, such as the heat conduction problem with a thin insulating layer or with a phase change at the material interface, and the percolation problem with a filter membrane, which also receive much attention. 
The solutions and fluxes of these problems are discontinuous at the material interface and they satisfy some jump conditions.
Such problems, which are also studied in Refs.~\cite{Wu,HOU}, can be formulated by the following equations:
\begin{equation}\label{jump_equ}
\begin{cases}
-\nabla\cdot\left(\kappa(x,y)\nabla u(x,y)\right)=Q(x,y),&~(x,y)\in~\Omega=(-1,1)\times(-1,1),\\
    [\![u(x,y)]\!]_{\Gamma}=\Phi(x,y),&~(x,y)\in~\Gamma,\\ 
     [\![\kappa(x,y)\nabla u(x,y)\cdot\boldsymbol{n}]\!]_{\Gamma}=\Psi(x,y),&~(x,y)\in~\Gamma,\\
    u(x,\pm1)=\ln(1+x^2),&~x\in [-1,1],\\
    u(\pm1,y)=\ln(1+y^2),&~y\in [-1,1],
\end{cases}
\end{equation}
where the coefficient $\kappa(x,y)$ and the source term $f(x,y)$ are as follows:
\begin{align}
    &\kappa(x,y)=
    \begin{cases}
        \cos(x+y)+2,~ & (x,y)\in\Omega_1=\{(x,y)|x^2+y^2<0.5^2\},\\
        \sin(x+y)+2,~ &(x,y)\in\Omega \backslash \Omega_1,
    \end{cases}\\
    &Q(x,y)=
    \begin{cases}
        4(\cos(x+y)+1)\sin(x+y),~ & (x,y)\in\Omega_1,\\
        -2\cos(x+y)\frac{x+y}{x^2+y^2},~ &(x,y)\in\Omega \backslash \Omega_1.
    \end{cases}
\end{align}
In Eq.~\eqref{jump_equ}, the interface $\Gamma$ is a circle with a radius of 0.5 and centered at $(0,0)$. 
Note that $[\![\mu]\!]_{\Gamma}:=\mu|_{\Gamma^+} - \mu|_{\Gamma^-}$ denotes the jump of $\mu$ across the interface. 
$\Phi(x,y)$ and $\Psi(x,y)$ can be derived from the exact solution below.

The exact solution of this case is 
\begin{equation}
    u(x,y)=
    \begin{cases}
        \sin(x+y),~ & (x,y)\in\Omega_1,\\
        \ln(x^2+y^2),~ &(x,y)\in\Omega \backslash \Omega_1.
    \end{cases}
\end{equation}

\begin{figure}[!htbp]
\centering
\includegraphics[width=0.9\textwidth]{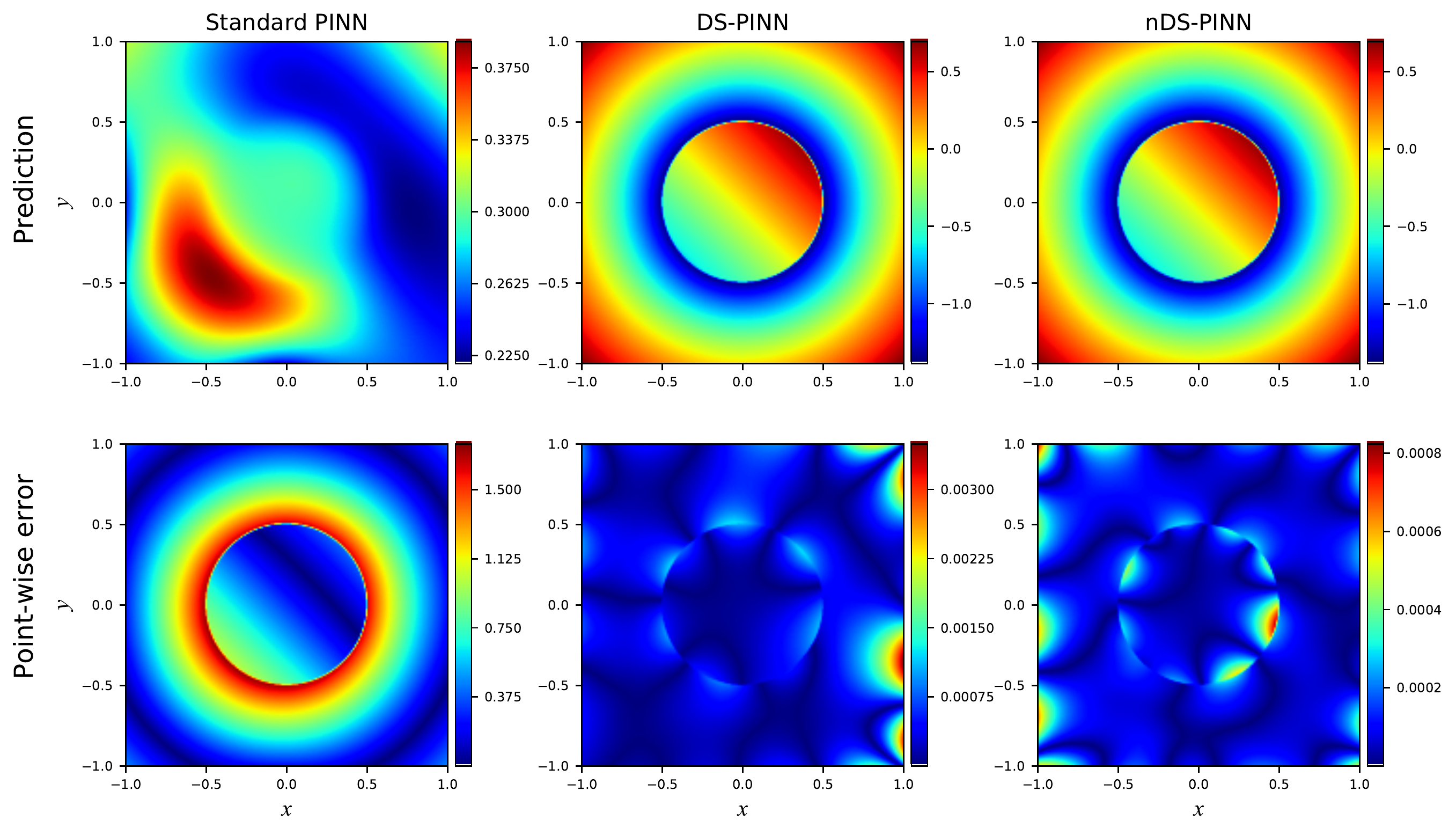}
\caption{Predictions and point-wise errors of three PINN methods for Sect.~\ref{Ex5}. }
\label{jumpresultspic}
\end{figure}

For this model, our methods are fully applicable, with only a slight modification of the loss term $\mathcal{L}_{\Gamma}(\theta;\tau_{\Gamma})$ 
by replacing the continuity conditions with the jump conditions.
The computational results using different PINN methods are shown in Figure \ref{jumpresultspic} and Table \ref{jumpresulttab}. It can be seen that the standard PINN method is powerless for such a model, while our methods DS-PINN and nDS-PINN (with separation distance ${d}=2$) can solve this model exactly, especially the nDS-PINN method gives extremely accurate computational results.

It should be emphasized that for the nDS-PINN method, since $\Phi(x,y)$ and $\Psi(x,y)$ in the jump condition are known functions, we used them to normalize the interface loss term $\mathcal{L}_{\Gamma}(\theta;\tau_{\Gamma})$.
The result of our nDS-PINN method is consistent with that of the INN method, which is taken from Ref.~\cite{Wu}.

\begin{table}[!htbp]
		\setlength{\abovecaptionskip}{0cm}
		\setlength{\belowcaptionskip}{0.2cm}
	\centering
	\caption{$\left\|e\right\|_{\mathbb{L}_2}$ errors of four PINN methods for Sect.~\ref{Ex5}.}
	\begin{tabular}{c|c}
		\hline
		Method  &  $\left\|e\right\|_{\mathbb{L}_2}$ \\ 
		 \hline
		 Standard PINN        &   $1.33\pm0.02\times10^{0}$\\
         INN~\cite{Wu}      &    $5.6\times10^{-4}$\\
		 DS-PINN (this work)      &    $1.7\pm0.8\times10^{-3}$\\
		 nDS-PINN (this work)    &   $3.8\pm1.2\times10^{-4}$\\
		 \hline
	\end{tabular}
 \label{jumpresulttab}
\end{table}

\section{Conclusions}\label{Results}
For a class of multi-material diffusion problems, this paper first analyzed the reasons why the standard PINN cannot be applied; then derived two continuity conditions that should be satisfied at the material interface, the use of which can effectively fill in the missing information at the interface; further, we designed a domain separation strategy to overcome the problem that the solution function cannot be expressed by a single neural network due to the discontinuity of its derivatives at the interface. Finally, by combining the above  two works, we improved the standard PINN by adding special terms to the loss function so that the interface conditions are accurately represented in a single neural network, which makes the obtained prediction function fully reflect the characteristics of the solution at the interface, giving very accurate predictions near the interface. 
In addition, we design a problem-adapted normalization method for the loss term, which can further significantly improve the accuracy of the prediction.
Various numerical experiments verify the effectiveness of our method. The new method perfectly solves the problem that the standard PINN cannot be adapted to the multi-material diffusion model.
We believe that this work provides a novel idea for PINN to solve partial differential equations with non-smooth solutions, and it is a useful development of the standard PINN.

Note that the methods in this paper are only for linear multi-material diffusion equations. It is our future work to study PINN methods for solving nonlinear multi-material diffusion problems.

\section*{Code availability}

The code of this work is publicly available online via  \url{https://doi.org/10.5281/zenodo.7927544}.

\section*{Acknowledgements}
The work is supported by the National Science Foundation of China under Grant No.12271055, the Foundation of CAEP (CX20210044), the Natural Science Foundation of Shandong Province No.ZR2021MA092, and the Foundation of Computational Physics Laboratory.

\bibliography{mybibfile}

\end{document}